\documentclass[12pt]{article}
\usepackage{euscript}
\usepackage{a4wide}
\usepackage[T2A]{fontenc}
\usepackage[utf8]{inputenc}
\usepackage{amsfonts}
\usepackage{amssymb, amsthm}
\usepackage{amsmath}
\usepackage{mathtools}
\usepackage{needspace}
\usepackage{lipsum}
\usepackage{comment}
\usepackage{cmap}
\usepackage[pdftex]{graphicx}
\usepackage{hyperref}
\usepackage{epstopdf}
\usepackage[matrix, arrow, 	curve]{xy}
\usepackage{amscd}
\usepackage{esint}

\begin{document}

\renewcommand{\proofname}{Proof}

\renewcommand{\d}{\partial} 
\newcommand{\Z}{\mathbb{Z}}
\newcommand{\N}{\mathbb{N}}
\newcommand{\R}{\mathbb{R}}
\newcommand{\Q}{\mathbb{Q}}
\newcommand{\K}{\mathbb{K}}
\newcommand{\Cm}{\mathbb{C}}
\newcommand{\Pm}{\mathbb{P}}
\newcommand{\B}{\mathcal{B}}
\newcommand{\Zero}{\mathbb{O}}
\newcommand{\ilim}{\int\limits}
\newcommand{\slim}{\sum\limits}
\newcommand{\action}{\curvearrowright}
\newcommand{\E}{\mathbb{E}}
\newcommand{\BB}{\overline{B}}
\newcommand{\D}{\mathcal{D}}
\newcommand{\T}{\mathbb{T}}
\newcommand{\F}{\mathcal{F}}
\newcommand{\Sf}{\mathbb{S}}
\newcommand{\Vol}{\mbox{V}}
\newcommand{\mint}{\strokedint\limits}
\newcommand{\const}{\mbox{const}}
\newcommand{\supp}{\mbox{supp}}
\newcommand{\dist}{\mbox{dist}}
\newcommand{\Hess}{\mbox{Hess}}
\newcommand{\Ker}{\mbox{Ker}}
\newcommand{\Hd}{\mathcal{H}}
\renewcommand{\div}{\mbox{div}}
\newcommand{\diam}{\mbox{diam}}
\newcommand{\NTlim}{\mbox{n.t.lim }}
\newcommand{\mydet}{\mbox{det}}
\newcommand{\Id}{\mbox{Id}}

\theoremstyle{plain}
\newtheorem{thm}{Theorem}
\newtheorem{lm}{Lemma}
\newtheorem*{st}{Statement}
\newtheorem*{prop}{Properties}
\newtheorem*{cl}{Claim}

\theoremstyle{definition}
\newtheorem{defn}{Definition}
\newtheorem{ex}{Ex}
\newtheorem{cor}{Corollary}

\theoremstyle{remark}
\newtheorem{rem}{Rem}
\newtheorem*{note}{Note}

\title{On an obstacle to the converse of Dahlberg's theorem in high codimensions}
\author{Polina Perstneva \thanks{The author was partially supported by the Simons Foundation grant 601941, GD.} \\ 
}
\date{\today}

\maketitle

\begin{abstract}
    It has been recently understood that the harmonic measure on the boundary $E = \partial \Omega$ of a domain $\Omega$ in $\R^n$ is absolutely continuous with respect to the Hausdorff measure $\Hd^{n - 1}$ on $E$ if and only if the boundary $E$ is rectifiable. Then, by G. David, M. Engelstein, J. Feneuil, S. Mayboroda and other coauthors, a notion of harmonic measure for Ahlfors-regular sets $E$ of higher codimension $n - d$ was developed with the aid of the operator $L_\alpha = -\div D_{\alpha}^{-n + d + \alpha} \nabla$, where $\alpha > 0$ and $D_\alpha$ is a certain regularized distance function to the set $E$. A program was launched to establish analogous to the classical case equivalence between rectifiability of the higher-codimensional set $E$ and good relations of the (new) harmonic and Hausdorff measures. The sufficiency of rectifiability for quantitative absolute continuity was only just obtained. For the other direction the main obstacle is to prove that, roughly, the equation $L_\alpha D_\alpha = 0$ is true only when the set $E$ is a hyperplane. In this paper we prove some first results which indicate that the latter conjecture may be true. We also explain that a certain natural strategy to tackle the problem does not work till the end.
\end{abstract}

\tableofcontents

\section{Introduction}


This note 
belongs to a long tradition of studying the relations between the geometry of a domain $\Omega \subset \R^n$
and the analysis on $\Omega$ or $E = \d \Omega$, in particular relative to second order elliptic operators on $\Omega$.
Important examples of this are the study of the absolute continuity of the harmonic measure on $E$ with respect to the surface measure.
This subject has a long history (and we refer the reader to the introduction of \cite{GSJ18} for a thorough survey); let us just
mention here two important results in the spirit of our interests here. In \cite{Da}, B. Dahlberg showed that when $\Omega$ is a Lipschitz domain, the harmonic measure on $E$ is mutually absolutely continuous with respect to the surface measure, and even given by an $A_\infty$ weight; there were lots of important results before, but mostly when $n=2$ and related to conformal mappings.
After this, it was slowly understood that the main issue in this problem was the rectifiability of the boundary $E$. Recall that $d$-rectifiability of a set means that it can be represented as at most countable union of Lipschitz graphs (of dimension $d$) and a set of zero Hausdorff measure $\Hd^d$. Later, the technology
improved to the point that one could also worry about the converse results, i.e., what can be said about $E$ when
the harmonic measure is absolutely continuous. In 2015, after a long series of works, the question was finally settled
in \cite{AHM3TV} by J. Azzam, S. Hoffman, M. Mourgoglou, J. M. Martell, S. Mayboroda, X. Tolsa and A. Volberg. They proved
in particular that for $n \geq 2$ and a set $E$ with $\Hd^{n-1}(E) < \infty$, the absolute continuity of the harmonic measure on $E$ with respect to $\Hd^{n - 1}$ on $E$ is equivalent to its rectifiability. 

Many of the results above, which initially concerned the Laplacian $\Delta$, were also extended to a class of elliptic operators $L$ that are sufficiently close to constant coefficient elliptic operators. See for instance  $\cite{AGMT}$, $\cite{HMM2}$.


After these successes, G. David, J. Feneuil, and S. Mayboroda \cite{GSJ18} 
started to inquire if the same philosophy is still true for domains $\Omega$ with a lower-dimensional boundary, and more precisely domains $\Omega \subset \R^n$ 
such that $E = \d \Omega$ is Ahlfors regular with a dimension $d < n-1$. Recall that $E$ is called $d$-Ahlfors regular if for some $C_0 \ge 1$ the double inequality 
$$C_0^{-1}r^d \le \Hd^d(E \cap B(x, r)) \le C_0 r^d$$
is true whenever $x \in E$ and $0 < r < \diam(E)$. Notice that in such a case there is no complementary component, i.e., $\Omega = \R^n \setminus E$,
and it is checked in \cite{GSJ18} that $\Omega$ has non-tangential access. 
But the harmonic measure (associated to $\Delta$) is not well defined
on $E$, for instance because Brownian paths do not see $E$, so the authors had to use a different class of (degenerate) elliptic
operators, adapted to the geometry of $E$. They consider divergence form operators
\begin{equation}\label{frst}
L = - \div A \nabla,
\end{equation}
where the matrix valued function $A: \Omega \to M_n(\R)$ satisfies the modified ellipticity conditions
$$\delta(x)^{n - d - 1}A(x) \zeta \cdot \xi \le C_1 |\zeta||\xi|, x \in \Omega, \zeta, \xi \in \R^n, \quad \mbox{and}$$
\begin{equation}\label{scnd}
\delta(x)^{n - d - 1}A(x)\zeta \cdot \zeta \ge C_1^{-1}|\zeta|^2, x \in \Omega, \zeta \in \R^n,
\end{equation}
where $C_1 \ge 1$ is a constant, and $\delta(x) = \dist(x, E)$ is the distance function from $x \in \Omega$ to $E$.  
With these operators $L$, \cite{GSJ18} establishes an analogue of the usual theory of elliptic operators: the existence and uniqueness of solutions, and some regularity of the latter. This allows one to define a (degenerate) elliptic measure on the boundary set $E$ in a ``usual'' way, along with a notion of a Green function. 


Then they worry about analogues in this higher co-dimension context of the direct absolute continuity results, like 
Dahlberg's theorem above. The concerned operators have to lie in a much smaller class; recall that in the classical case, 
typical absolute continuity results concern only perturbations of the Laplacian. They settle on the nicest operator they found,
namely
\begin{equation}\label{3}
L = L_{\alpha,\mu} = - \div D_{\alpha,\mu}^{-n + d + 1} \nabla,
\end{equation}
where $\mu$ is some $d$-dimensional Ahlfors regular measure on $E$, $\alpha > 0$ is a parameter, and
the corresponding smooth distance function $D_{\alpha,\mu}$ is defined by
\begin{equation}\label{mu_dist}
    D_{\alpha, \mu}(x) = \left(\ilim_E{|x - y|^{-d - \alpha} d\mu(y)}\right)^{-1/\alpha}.
\end{equation}
We say that the measure $\mu$ is a $d$-dimensional Ahlfors regular measure on $E$ when the (closed) support of
$\mu$ is $E$, and when there is a constant $C_0$ such that 
\begin{equation} \label{?5}
C_0^{-1}r^d \le \mu(B(x, r)) \le C_0 r^d
\end{equation}
for $x \in E$ and $0 < r < \diam(E)$. It is easy to check that since $\mu$ is Ahlfors regular, $D_{\alpha,\mu}(x)$ is equivalent to
$\dist(x, E)$, hence $ L_{\alpha,\mu}$ satisfies the constraints \eqref{frst} and \eqref{scnd}.

It is well known that when there is a $d$-Ahlfors regular measure on $E$, then the restriction $\Hd^d_{\vert E}$ 
of the Hausdorff measure is $d$-Ahlfors regular too. In fact, in \cite{GSJ17}, the authors restrict to $\mu = \Hd^d_{\vert E}$, and prove that when $E$ is the graph of a Lipschitz function with a small enough Lipschitz norm, the harmonic measure
associated to $L_{\alpha,\mu}$ is mutually absolutely continuous with respect to $\mu$, and given by an $A_\infty$
weight. The result extends to any $d$-Ahlfors regular measure $\mu$ on $E$, and later on it was proved in 
\cite{GS20} and \cite{J20} that this result extends to the case where $E$ is uniformly rectifiable of dimension $d < n - 1$, $\mu$ is any
$d$-Ahlfors regular measure on $E$, and $\alpha > 0$. 


At this point it makes sense to look for a converse, but it was found in \cite{GMS20} that the following anomaly occurs. For the ``magic'' number $\alpha = n - d - 2$, the function $\eqref{mu_dist}$ is a solution to $L_{\alpha, \mu} \cdot = 0$ on $\Omega$ for the operator $\eqref{3}$. This implies that the harmonic measure for this operator is absolutely continuous with respect to Hausdorff measure $\Hd^d$ on $E$, no matter how irregular geometrically this set really is. This means that the full analogue of the reverse to Dahlberg's theorem in higher codimension cannot be true. Yet it sounds reasonable to expect that, except in the magic case when $d < n - 2$ and $\alpha = n - d + 2$,
the uniform rectifiability of $E$ follows from the $A_\infty$-absolute continuity of the harmonic measure associated to $L_{\mu,\alpha}$
with respect to $\mu$.

In \cite{GS_f20}, the authors propose to address the different issue of good approximation of the Green function for $L_{\mu,\alpha}$
(with a pole at $\infty$) by distance functions to $E$. They prove some direct results,
and also show that some interesting 
converse results will follow if one proves the following conjecture. 
Let $E$ be a $d$-Ahlfors regular set, $D_{\alpha, \mu}$ -- the regularized distance function as above in $\eqref{mu_dist}$, and $L_{\alpha, \mu}$ -- the degenerate elliptic operator mentioned before. Then 
\begin{equation}\label{solution}
    L_{\alpha, \mu} D_{\alpha, \mu} = 0 \quad \mbox{in} \quad \Omega = \R^n \setminus E 
\end{equation}
is never true except for the following two cases:
\begin{enumerate}
    \item when $d < n - 2$ and $\alpha = n - d - 2$, 
    \item when $E = \R^d$ for some integer $d$ and $\mu = c\Hd^d|_E$ for some positive constant $c$.
\end{enumerate}

In this work we will make the first step in the study of this conjecture. Before we state our results, let us explain our motivation and give some definitions. 

Our global goal is to prove that, for the case when $\alpha$ is not the ``magic'' number, the only possible solution to $\eqref{solution}$ is what we call the flat solution. We wanted to start with explaining why there are no solutions in a neighbourhood of the flat one. This could be easier than the study of global solutions, since we can view the non-flat ones as small perturbations of the flat. Essentially, any story about perturbations involves a parameter and a family of solutions corresponding to it. This is why an often-used plan to prove the absence of solutions in a neighbourhood is to prove first the absence of parametric families of solutions. Which is exactly what we will do. Then the transition to the absence of individual solutions is usually not too hard, though in our situation the most logical scheme does not work: we will discuss this in section 5. So the results we state and prove below are actually the best one could do trying to follow the described plan to solve the hypothesis.

Let us discuss in more details what we mean by parametric families in a neighbourhood of the flat solution. To start with, consider the case when the measures $\mu$ of all of the functions $D_{\alpha, \mu}$ of our family live on the hyperplane $E = \R^d$, but their densities with respect to the Hausdorff measure are not constants. Then the easiest one-parameter family of solutions of $\eqref{solution}$ to study is $\{D_{\alpha, \mu_t}\}_{t \in [0, t_0)}$ with $E = \R^d$ and $\mu_t = (1 + t\phi)d\Hd^d$ for a fixed function $\phi \in L^\infty(\R^d)$. The next step is to switch to the solutions $D_{\alpha, \mu}$ such that $\mu$ is supported on graphs close to the hyperplane $\R^d$. Here the easiest one-parameter family to study is $\{D_{\alpha, \mu_t}\}_{t \in [0, t_0)}$ with $E_t = Im(Id + t\psi)$, and $(Id + t\psi)^{-1}(\mu_t) = (1 + t\phi)d\Hd^d$ for fixed functions $\phi \in L^\infty(\R^d) \cap L^1(\R^d)$, $\psi \in Lip(\R^d) \cap L^1(\R^d)$. These are the examples one can keep in mind. 

Our results concern more general families of solutions. We start with measures supported on $E = \R^d$.

\begin{defn}\label{defn1}
We call a non-trivial one-parameter differentiable family of flat perturbations of the flat solution a family of solutions $\{D_{\alpha, \mu_t}\}_{t \in [0, t_0)}, t_0 > 0$, of $\eqref{solution}$ such that for any $t$ the measure $\mu_t$ is supported on $E$ has density $1 + \phi_t, \phi_t \in L^\infty(\R^d)$ with respect to $\Hd^d$ with the following properties:
\begin{enumerate}
    \item $\phi_0 = 0$,
    \item the family of densities is Frechet differentiable at zero in $BMO (\R^d)$; that is, there exists a function $\frac{\partial \phi_t(y)}{\partial t}|_{t = 0} = \frac{\partial \phi_t}{\partial t}(0, y)$ such that
    $$\|\phi_t - t \frac{\partial \phi_t}{\partial t}(0, \cdot)\|_{BMO} = o(t), \quad \mbox{as} \; t \to 0,$$
    \item the derivative $\frac{\partial \phi_t}{\partial t}(0, \cdot)$ is a non-constant function,
    \item $$\|\phi_t\|_{BMO}, \Big|\ilim_{B(0, 1)}{\phi_t(y)dy}\Big| \le Ct \quad \mbox{and} \quad \Big| \ilim_{B(0, 1)}{\left(\phi_t(y) - t\frac{\partial \phi_t}{\partial t}(0, y)\right)dy} \Big| = o(t), \; t \to 0.$$
\end{enumerate}
\end{defn}
The integral conditions above on the ball $B(0, 1)$ help, because the $BMO$ norm itself does not control averages on large balls. In the next definition, when we write that the modulus of a vector-valued function or its norm (in $BMO$) admits an estimate, we mean that the moduli or norms of every component of this vector admit it.

\begin{defn}\label{defn2}
We call a non-trivial one-parameter differentiable family of graph perturbations of the flat solution a family of solutions $\{D_{\alpha, \mu_t}\}_{t \in [0, t_0)}, t_0 > 0$, of $\eqref{solution}$ such that for any $t$ the support $E_t$ of $\mu_t$ is the image of a Lipschitz function $Id + \psi_t: \R^d \to \R^n$, where $\psi_t: \R^d \to \R^{n - d}$ has the Lipschitz constant $Ct$, the measure $\mu_t$ is the image by $Id + \psi_t$ of $(1 + \phi_t)d\Hd^d, \phi_t \in L^\infty(\R^d)$, and the following conditions hold:
\begin{enumerate}
    \item $\phi_0 = \psi_0 = 0$,
    \item the families of densities $\{\phi_t\}$ and graph functions $\{\psi_t\}$ are Frechet differentiable at zero in $BMO(\R^d)$: there exist functions  $\frac{\partial \phi_t(y)}{\partial t}|_{t = 0} = \frac{\partial \phi_t}{\partial t}(0, y)$ and $\frac{\partial \psi_t(y)}{\partial t}|_{t = 0} = \frac{\partial \psi}{\partial t}(0, y)$ in $BMO(\R^d) \cap L^1(\R^d)$ such that
    $$\|\phi_t - t\frac{\partial \phi_t}{\partial t}(0, \cdot)\|_{BMO} = o(t), \quad \mbox{as} \; t \to 0,$$
    $$\|\psi_t - t\frac{\partial \psi_t}{\partial t}(0, \cdot)\|_{BMO} = o(t), \quad \mbox{as} \; t \to 0,$$
    \item the derivative $\frac{\partial \psi}{\partial t}(0, y)$ is not constant,
    \item for $F_t = \phi_t$ and $\psi_t$ 
    $$\|F_t\|_{BMO}, \Big|\ilim_{B(0, 1)}{F_t(y)dy}\Big| \le Ct \quad \mbox{and} \;  \Big| \ilim_{B(0, 1)}{\left(F_t(y) - t\frac{\partial F_t}{\partial t}(0, y)\right)dy} \Big| = o(t), \; t \to 0.$$
\end{enumerate}
\end{defn}

We think that the result we state below for the non-trivial one-parameter differentiable families of graph perturbations is true without the assumption that the derivatives $\frac{\partial \phi_t}{\partial t}(0, y)$ and $\frac{\partial \psi_t}{\partial t}(0, y)$ are in $L^1(\R^d)$ (or without any other similar summability assumption), but we did not manage yet to think of a better argument than the one in Subsection 5.2, which uses the Fourier transform.

We are now ready to state our first two main theorems. Recall that we are interested in the case when for our parameters $n, d$ and $\alpha$ we have $n - d > 2$ and $\alpha > 0$ is not ``magic'' ($\alpha \neq n - d - 2$), and we do not pose any additional restrictions on them in Theorems $\ref{mainthm_flat}$ and $\ref{mainthm_graph}$.  
 
\begin{thm}\label{mainthm_flat}
For and integer $d < n - 2$ and $E = \R^d$ there are no non-trivial one-parameter differentiable families of flat perturbations of the solution $D_{\alpha, \mu}$ with $\mu = c\Hd^d|_E$ of the equation $L_{\alpha, \mu} D_{\alpha, \mu} = 0$.
\end{thm}

\begin{thm}\label{mainthm_graph}
For any integer $d < n - 2$ there are no non-trivial one-parameter differentiable families of graph perturbations of the solution $D_{\alpha, \mu}$ with $\mu = c\Hd^d|_E$ and $E = \R^d$ of the equation $L_{\alpha, \mu} D_{\alpha, \mu} = 0$. 
\end{thm}

For the case when $E$ is a hyperplane we are able to provide another result in the spirit of non-existence of global solutions. We show that, if the density of the measure $\mu$ in $D_{\alpha, \mu}$ with respect to the Hausdorff measure is regular enough and it is not a constant, then $\eqref{solution}$ cannot be true. To make precise the notion of ``regular enough'' we remind the reader yet another definition.

\begin{defn}
We say that a function $f$ on $\R^d$ is in H\"older class $C^{k, \gamma}(\R^d)$ if it has continuous derivatives up to the order $k$ and the $k$th partial derivatives are H\"older continuous with the exponent $\gamma$, $0 < \gamma < 1$. A function $g$ is H\"older continuous with the exponent $\gamma$ if 
$$\sup_{x \neq y}{\frac{|g(x) - g(y)|}{|x - y|^\gamma}} < \infty.$$
\end{defn}

Additional restrictions on $\alpha$ and $d$ are posed due to our method. They come from the assumption of integrability of certain functions.

\begin{thm}\label{mainthm_direct}
If $E$ is a hyperplane of dimension $d$ ($E = \R^d$) such that $n - d > 4$, $\alpha > 2 + \varepsilon_0$ for some $0 < \varepsilon_0 < 1$, and the density of the measure $\mu$ with respect to the Hausdorff measure $\Hd^d$ on $E$ is not a constant, but of class $C^{2, \varepsilon}$ for some $0 < \varepsilon < \varepsilon_0$, then the function $D_{\alpha, \mu}$ as in $\eqref{mu_dist}$ cannot be a solution for the equation $\eqref{solution}$.
\end{thm}

The paper is organized as follows. In Section 2 we discuss a suitable reformulation of the equation $\eqref{solution}$ in terms of the Laplacian, and then give a representation of any harmonic function on $\Omega$ with certain asymptotics. In Section 3 we discuss non-tangential limits of the smooth distance function $D_{\alpha, \mu}$. In Section 4 we prove some facts we need about the space of functions with bounded mean oscillations. In Section 5 we prove Theorems $\ref{mainthm_flat}$ and $\ref{mainthm_graph}$, and discuss why our method does not seem to allow one to finish the proof of the conjecture about the uniqueness of flat solutions of $\eqref{solution}$. In Section 6 we prove Theorem $\ref{mainthm_direct}$.

~\

\subsubsection*{Acknowledgements} I am very grateful to my PhD thesis advisor, Professor Guy David, for the patient guidance and permanent support during the development of this paper. I also would like to thank Pierre-Gilles Lemari\'e-Rieusset and Ioann Vasilyev for fruitful discussions.

\section{Harmonic functions outside a $d$-Ahlfors regular set $E$}

\subsection{An observation}

Our method is based on a simple observation, which we now present. Throughout the text we suppose that $E$ is a $d$-Ahlfors regular set and $\mu$ is a $d$-Ahlfors regular measure on $E$. Suppose that the function $D_{\alpha, \mu}$ in $\eqref{mu_dist}$, which is easily seen to be smooth in $\Omega = \R^n \setminus E$,
is a solution to the equation $L_{\alpha, \mu} \cdot = -\div D_{\alpha, \mu}^{n - d - 1} \nabla \cdot = 0$ on $\Omega$. Then
\begin{equation}\label{1}
    0 = -L_{\alpha, \mu} D_{\alpha, \mu} = \div\left(D_{\alpha, \mu}^{-n + d + 1} \nabla D_{\alpha, \mu}\right) = D_{\alpha, \mu}^{- n + d + 1}\Delta D_{\alpha, \mu} + (- n + d + 1)D_{\alpha, \mu}^{- n + d}\slim_{i = 1}^n{\left(\frac{\partial}{\partial x_i}D_{\alpha, \mu}\right)^2}.
\end{equation}
For $\gamma \in \R$, we can compute 
\begin{equation}\label{2}
    \Delta D_{\alpha, \mu}^\gamma = \slim_{i = 1}^n{\left(\gamma D_{\alpha, \mu}^{\gamma - 1}\frac{\partial^2}{\partial x_i^2}D_{\alpha, \mu} + \gamma(\gamma - 1)D_{\alpha, \mu}^{\gamma - 2}\left(\frac{\partial}{\partial x_i}D_{\alpha, \mu}\right)^2\right)} \quad \mbox{on} \; \Omega.
\end{equation}
Evidently, if we pick $\gamma$ equal to $-n + d + 2$ the right-hand sides of $\eqref{1}$ and $\eqref{2}$ will coincide. Therefore $D_{\alpha, \mu}$ is a solution of $\eqref{solution}$ if and only if $D_{\alpha, \mu}^\gamma$ is harmonic outside $E$ for $\gamma = - n + d + 2$:
\begin{equation}\label{maineq}
    L_{\alpha, \mu}D_{\alpha, \mu} = 0 \quad \Longleftrightarrow \quad \Delta D_{\alpha, \mu}^\gamma = 0.
\end{equation}
This fact will provide us with an alternative representation for the integral $\eqref{mu_dist}$. We will explain this after providing the necessary preliminaries.

Note that for the ``magic'' number $\alpha = n - d - 2$ it is always true that $\Delta D_{\alpha, \mu}^\gamma = 0$ for the chosen exponent $\gamma$.

\subsection{The Newton potential}

Before we proceed, recall that our dimension $d$ is such that $n - d - 2 > 0$, and that we denote by $\delta(x)$ the distance function $\dist(x, E)$. Let $\sigma$ be the measure $\Hd^d$ restricted to $E$. The purpose of this subsection is to prove the following theorem.

\begin{thm}\label{ker_repr}
Let $u$ be a function harmonic outside the $d$-Ahlfors regular set $E \subset \R^n$ such that $|u(x)| \le C\delta(x)^{-n + d + 2}$. Then there exists a function $f \in L^\infty(E)$ such that 
$$u(x) = \ilim_{\R^n}{\frac{f(y) d\sigma(y)}{|x - y|^{n - 2}}}.$$
\end{thm}

This, of course, looks like the well-know representation of a solution of the equation $\Delta \cdot = 0$ as a convolution with the fundamental solution. But the latter is almost always used in the situation when the set $E$ is compact and has codimension one, while here the issue really is to prove that the Laplacian of $u$, in the sense of distributions, is an upper $d$-Ahlfors regular measure.

The next lemma uses standard arguments, but we give the proof for the sake of completeness. Following the traditions of PDE texts, from now on we usually denote all the various constants by the letter $C$.

\begin{lm}\label{newton}
Let $f$ be a function in $L^\infty(E, d\sigma)$. Consider 
$$u_f(x) = \ilim_{\R^n}{\frac{f(y)}{|x - y|^{n - 2}}d\sigma(y)};$$
with the assumptions above, the following holds:
\begin{enumerate}
\item $u_f$ is locally integrable ($u_f \in L_{1, loc}(\R^n)$),
\item $u_f$ is harmonic in $\R^n \setminus E$,
\item if $E$ is compact, then for $x$ such that $dist(x, E) = \delta(x) < \diam(E)$, and if $E$ is non-compact, everywhere on $\R^n \setminus E$
$$|u_f(x)| \le c_1\delta(x)^{- n + d + 2},$$
\item in the sense of distributions,
$$\Delta u_f = - |S_1|(n - 2)f d\sigma,$$
where $|S_1|$ is the area of the sphere of radius $1$ in $\R^n$.
\end{enumerate}
\end{lm}

\begin{note}
The potential $\ilim_{\R^n}{\frac{f(y)}{|x - y|^{n - 2}}d\sigma(y)}$ is often called the Newton potential. 
\end{note}

\begin{proof}
\begin{enumerate}
\item For all $r > 0$
$$\ilim_{|x| < r}{\ilim_{\R^n}{\frac{|f(y)| d\sigma(y)}{|x - y|^{n - 2}}}dx} \le \|f\|_\infty \ilim_{|x| < r}{\left(\ilim_{|y| < 2r}{\frac{d\sigma(y)}{|x - y|^{n - 2}}} + \ilim_{|y| \ge 2r}{\frac{d\sigma(y)}{|x - y|^{n - 2}}}\right)dx}.$$
The first term we estimate the following way:
$$\ilim_{|x| < r}{\ilim_{|y| < 2r}{\frac{d\sigma(y)}{|x - y|^{n - 2}}}dx} \le  \ilim_{|y| < 2r}{\ilim_{|x| < r}{\frac{dx}{|x - y|^{n - 2}}}d|\sigma|(y)} \le \ilim_{|y| < 2r}{\ilim_{|z| < r + 2r}{\frac{dz}{|z|^{n - 2}}}d\sigma(y)} $$ 
$$ = \sigma(E \cap B(0, r))\ilim_{|z| < 3r}{\frac{dz}{|z|^{n - 2}}} = \sigma(E \cap B(0, r))|S_1|\ilim_0^{3r}{\rho d\rho} = C r^d |S_1|\frac{(3r)^2}{2}.$$
To estimate the second term we exploit the fact that $x$ and $y$ are now far away from each other, and we also split the space into layers like always:

$$\ilim_{|x| < r}{\ilim_{|y| \ge 2r}{\frac{d\sigma(y)}{|x - y|^{n - 2}}}dx} \le \ilim_{|x| < r}{\slim_{k = 1}^\infty{\ilim_{2^k r \le |y| < 2^{k + 1}r}{\frac{d\sigma(y)}{|x - y|^{n - 2}}}}dx} $$ 
$$ \le \ilim_{|x| < r}{\slim_{k = 1}^\infty{(2^k r)^{- n + 2}C(2^{k + 1}r)^d} dx} \le C r^{d + 2},$$
since $- n + d + 2 < 0$.

\item It suffices to say that the function $\frac{1}{|x|^{n - 2}}$ is a fundamental solution of the equation $\Delta u = 0$ in $\R^n \setminus \{0\}$, and that 
$$\Delta u_f(x) = \ilim_{\R^n}{\Delta \frac{1}{|x - y|^{n - 2}} f(y)d\sigma(y)}.$$

\item Our strategy will be quite similar to the one we used for the proof of the local integrability:
$$|u_f(x)| \le \|f\|_\infty \left(\ilim_{B(x, 2\delta(x))}{\frac{d\sigma(y)}{|x - y|^{n - 2}}} + \slim_{k \ge 1}{\ilim_{B(x, 2^{k + 1}\delta(x))\setminus B(x, 2^k\delta(x))}{\frac{d\sigma(y)}{|x - y|^{n - 2}}}}\right).$$
The first integral inside the brackets can be estimated as
$$\ilim_{B(x, 2\delta(x))}{\frac{d\sigma(y)}{|x - y|^{n - 2}}} \le (2\delta(x))^{-n + 2}\ilim_{B(x, 2\delta(x))}{d\sigma(y)} \le (2\delta(x))^{-n + 2}c^d\delta(x)^d \le C \delta^{-n + 2 + d}.$$
For the terms in the second part,
$$\ilim_{B(x, 2^{k + 1}\delta(x))\setminus B(x, 2^k\delta(x))}{\frac{d\sigma(y)}{|x - y|^{n - 2}}} \le \left(2^{k}\delta(x)\right)^{-n + 2}(c2^{k + 1}\delta(x))^d = C (2^{-n + 2 + d})^k \delta^{-n + 2 + d}.$$
Thus, as we sum over $k > 0$, we get that the second term is equal to $C \delta^{-n + 2 + d}$, since $2^{-n + 2 + d} < 1$.

\item For any $\phi \in C_0^\infty$ we have
$$\langle \Delta u_f, \phi \rangle = \langle u_f, \Delta\phi \rangle$$
by the definition of the distribution $\Delta u_f$. But
$$\langle u_f, \Delta\phi \rangle = \ilim_{\R^n}{u_f(x)\Delta\phi(x)dx} = \ilim_{\R^n}{\ilim_{\R^n}{\frac{f(y) d\sigma(y)}{|x - y|^{n - 2}}}\Delta\phi(x)dx}.$$
Since we know that $u_f \in L_{1, loc}(\R^n)$, $\ilim_{\R^n}{|u_f(x)\Delta\phi(x)|dx} < \infty$, and by Fubini's theorem 
$$\ilim_{\R^n}{\ilim_{\R^n}{\frac{f(y) d\sigma(y)}{|x - y|^{n - 2}}}\Delta\phi(x)dx} = \ilim_{\R^n}{\ilim_{\R^n}{\frac{\Delta \phi(x)}{|x - y|^{n - 2}}dx}f(y) d\sigma(y)}.$$ 

Now we use the fact that the solution of the distributional equation $\Delta u = \delta_0$ is $-\frac{1}{(n - 2)|S_1|}\frac{1}{|x|^{n - 2}}$. This is well-known, but one could also check out $\cite{AG}$, p.103, Lemma 4.3.6, for example. Therefore
$$\phi(y) = \delta_0 * \phi(y) = \Delta \left(- \frac{1}{(n - 2)|S_1|}\frac{1}{|x^{n - 2}|}\right) * \phi(y) = -\frac{1}{|S_1|(n - 2)} \ilim_{\R^n}{\frac{\Delta \phi(x)dx}{|y - x|^{n - 2}}}.$$
We can integrate the expression above and conclude that
$$\ilim_{\R^n}{\ilim_{\R^n}{\frac{\Delta \phi(x)}{|x - y|^{n - 2}}dx} f(y) d\sigma(y)} = -(n - 2)|S_1|\ilim_{\R^n}{\phi(y) f(y) d\sigma(y)}.$$
Thus,
$$\langle \Delta u_f, \phi\rangle = -|S_1|(n - 2)\langle f d\sigma, \phi\rangle.$$
\end{enumerate}
\end{proof}

~\

Now we study an arbitrary function $u$ harmonic in $\R^n \setminus E$ and such that $|u(x)| \le C \delta(x)^{-n + 2 + d}$. For this we need to recall that $\mu$ is an upper $d$-Ahlfors regular measure on $E$ if there is a constant $C_0$ such that for $x \in E$ and $0 < r < \diam(E)$ one has $\mu(B(x, r)) \le C_0r^d$.

\begin{lm}\label{mes}
For a function $u$ as above, $u \in L_{1, loc}(\R^n)$ and the distribution $\Delta u$ is an upper $d$-Ahlfors regular measure on $E$.
\end{lm}
\begin{proof}
The first part, that $u$ lies in $L_{1, loc}(\R^n)$, can be viewed as a corollary of the following fact. For $\sigma = \Hd^d|_E$, one has $\delta(x)^{-n + d + 2} \le c u_1(x)$ for the function $u_1(x)$ as in lemma $\ref{newton}$. Indeed, let us pick a point $y_0$ such that $d(x, y_0) \le 2\delta(x)$ and $r = \delta(x)$. Then
$$u_1(x) = \ilim_{\R^n}{\frac{d\sigma(y)}{|x - y|^{n - 2}}} \ge \ilim_{E \cap B(y_0, r)}{\frac{d\sigma(y)}{|x - y|^{n - 2}}} \ge \ilim_{E \cap B(y_0, r)}{\frac{d \sigma(y)}{(|x - y_0| + |y_0 - y|)^{n - 2}}}$$ 
$$ \ge \ilim_{E \cap B(y_0, r)}{\frac{d\sigma(y)}{(3\delta(x))^{n - 2}}} = C \delta(x)^{- n + 2}\sigma(E \cap B(y_0, r)) \ge C\delta^{- n + 2}r^d = C\delta(x)^{- n + d + 2}.$$
Therefore for an arbitrary radius $r$, 
$$\ilim_{|x| < r}{|u(x)|dx} \le C \ilim_{|x| < r}{\delta(x)^{-n + 2 + d}dx} \le C\ilim_{|x| < r}{\ilim_{\R^n}{\frac{d\sigma(y)}{|x - y|^{n - 2}}}dx},$$
and now we can argue that in Lemma $\ref{newton}$ above we have already demonstrated that the function $u_1$ is locally integrable. One could also repeat the argument we used before, but for $\delta(x)^{-n + d + 2}$. 

\medskip

We now prove the second part of the statement. Observing that $-\Delta u$ is supported on $E$ is easy, since, if $\phi \in C_0^\infty(\R^n)$ has support outside $E$, then $\langle \Delta u, \phi\rangle$ is zero. For the upper-regularity we use an argument of approximate identity. Recall that, besides the local integrability of the function $u$, we know that its total mass inside a ball of small radius $r$ centered at $E$ is majorized up to a constant by $r^{d + 2}$, as we saw in Lemma $\ref{newton}$. Let $\{\phi_r\}$ be a standard approximate identity. We introduce the family of functions $\{u * \phi_r\}$. We will see that for each ball $B(x, \rho)$ centered at $E$ the total mass of the Laplacian of $u * \phi_r$ inside this ball is less than $C\rho^d$ if $r << \rho$, where $C$ does not depend on $r$. Therefore the same will be true for $\Delta u$ as a limit of $\Delta (u * \phi_r)$, which will conclude our proof. So, we fix $x \in E$, suppose that $r << \rho$ and estimate $\ilim_{B(x, \rho)}{|\Delta (u * \phi_r)|}$. Since $u * \phi_r$ is harmonic at distance at least $r$ from $E$, we really integrate over the set $B(x, \rho) \cap \{\delta(y) \le r\}$. We can cover this set by less than  $C\left(\frac{\rho}{r}\right)^d$ balls $B(x_i, r)$ centered at $E$. Indeed, let $\{B(x_i, r/5)\}_{i \in I}$ be any covering of $E \cap B(x, \rho)$. The Vitali lemma says that we can choose a finite $I_0 \subset I$ such that $B(x_i, r/5), i \in I_0$ do not intersect, and $\cup_{I_0} B(x_i, r)$ covers $E \cap B(x, \rho)$. The $d$-Ahlfors regularity of $E$ then implies that $|I_0|r^d \le C \rho^d$, which implies the bound on the number of balls in the covering $\{B(x_i, r)\}_{i \in I_0}$. The set $B(x, \rho) \cap \{\delta(y) < r\}$ can be covered by the $\{B(x_i, 2r)\}_{i \in I_0}$.

Using this development, and also that $\sup|\Delta \phi_r|$ is less than $r^{-2}$, we can estimate the total mass of the Laplacian of $u * \phi_r$ from above:
$$\ilim_{B(x, \rho)}{|\Delta (u * \phi_r)|} = \ilim_{B(x, \rho) \cap \{\delta \le r\}}{u * |\Delta \phi_r|} \le \slim_i{\ilim_{B(x_i, r)}{u * |\Delta \phi_r|}} \le \slim_{i}{\sup|\Delta \phi_r|\ilim_{B(x_i, 2r)}{u}}$$ 
$$ \le \slim_{i}{C r^{-2}r^{d + 2}} \le \left(\frac{\rho}{r}\right)^d Cr^d \le C\rho^d.$$

\end{proof}

\begin{lm}\label{mes2}
Let $\mu$ be a $d$-upper-regular measure on $E$, then it has a bounded density with respect to the measure $\sigma = \Hd^d|_E$. Moreover, if the measure $\mu$ is $d$-Ahlfors regular, then the density is also bounded away from zero.
\end{lm}
\begin{proof}
Let us check first that $\mu$ is absolutely continuous with respect to $\sigma$, or, equivalently, that if $\sigma(A) = 0$, then $\mu(A) = 0$ as well. Indeed, $\sigma(A) = 0$ if and only if for every $\varepsilon > 0$ exist a covering $\{B(x_i, r_i)\}$ of the set $A$ such that $\slim_{i}{r_i^d} < \varepsilon$. Then
$$|\mu(A)| \le |\mu(\cup B(x_i,r_i))| \le \slim_i{|\mu(B(x_i, r_i))|} \le C \slim_i{r_i^d} < C \varepsilon,$$
due to the upper-regularity of $\mu$. So, $\mu(A) \le C \varepsilon$ for every $\varepsilon > 0$, which implies that $\mu(A) = 0$. 

Now, it is well-known, see for instance \cite{M} Theorem 2.17, that the density function $f$ of the absolutely continuous part of a Radon measure $\mu$ with respect to a Radon measure $\sigma$ is (almost everywhere) equal to the function $D_{\sigma}\mu(x)$, the limit of $D_{\sigma}\mu(x, r) = \frac{\mu(B(x, r))}{\sigma(B(x, r))}$  with respect to $r \to 0$. In our case, since $\mu$ is absolutely continuous with respect to $\sigma$, the density function we are looking for is exactly $D_{\sigma}\mu$. It is bounded, since $\mu(B(x, r)) \le C r^d$ by the $d$-upper-regularity property. In the case then $\mu$ is $d$-Ahlfors regular the density function is also bounded away from zero, since $\mu(B(x, r)) \ge c r^d$.
\end{proof}

~\

\begin{proof}[Proof of theorem \ref{ker_repr}]
From lemma $\ref{mes}$ we know that if $u$ is harmonic outside $E$ and $|u(x)| \le C\delta(x)^{- n + d + 2}$, then $\Delta u$ is upper $d$-Ahlfors regular and therefore has a bounded density $\tilde{f}$ with respect to $\sigma$. Observe that for $f = \tilde{f}\frac{1}{|S_1|(n - 2)}$
$$\ilim_{\R^n}{(u - u_f) \Delta \phi} = \ilim_{\R^n}{\Delta (u - u_f) \phi} = - \ilim_{\R^n}{\phi \tilde{f} d\sigma} + \ilim_{\R^n}{\phi \tilde{f} d\sigma} = 0,$$
for every $\phi \in C_0^\infty$. Now the fundamental lemma of variational calculus tells us that $u$ and $u_f$ can differ only by a linear function. But if that linear function is not zero, this is impossible, since $|u(x)| \le C\delta(x)^{-n + d + 2}$, and this is arbitrarily small far away from $E$, which is not the case with $u = u_f + L$, $L$ linear and non-zero. 
\end{proof}

\section{Non-tangential limits}
The main goal of this section is to prove that, if $E$ is a $d$-Ahlfors regular and rectifiable set, and $u$ is a function harmonic outside $E$ and such that $|u(x)|$ is comparable to $\delta(x)^{-n + d + 2}$, then
the function $u(x)\delta(x)^{n - d - 2}$ has a non-tangential limit at almost every point $y_0 \in E$ equal to $c f(y_0),$ where $f$ is the density of the measure $\Delta u$ with respect to $\Hd^d|_E$, and the constant $c$ does not depend on $u$ or the point $y_0$. A statement very similar to this is also proved in $\cite{GMS20}$ (see section 5), and we use its techniques and follow mostly its exposition.  

Let us recall first the necessary terminology and definitions. With the assumption of rectifiability of the set $E$, at almost every with respect to the measure $\sigma = \Hd^d|_E$ point $y_0$ of $E$ we can find the unique tangent hyperplane $T_{y_0}E$ of dimension $d$: see, for example, $\cite{M}$ p. 219 Ex. 7 or $\cite{MShah}$ Ex. 41.21. Let $R > 0$ be small enough, and for the sets of points $y_0$ where $T_{y_0}E$ exists let us define for $\eta \in (0, 1)$ the set
$$\Gamma_{R, \eta} = \{x \in \left(\R^n \setminus E\right) \cap B(y_0, R): \; \dist(x, E) \ge \eta |x - y_0|\}.$$
We will call $\Gamma_{R, \eta}$ a non-tangential access region of the point $y_0$ (with an aperture $\eta$). Let $v$ be a function defined at least on $\R^n \setminus E$. The non-tangential limit of $v$ at a point $y_0 \in E$, denoted by $\NTlim_{x \to y_0} v(x)$, if it exists, is the common limit of $\{v(x_i)\}$, where $\{x_i\}$ is any sequence in $\Gamma_{R, \eta}$ such that $x_i$ tends to $y_0$ as $i$ tends to infinity.  

\begin{center}
\includegraphics[scale=0.7]{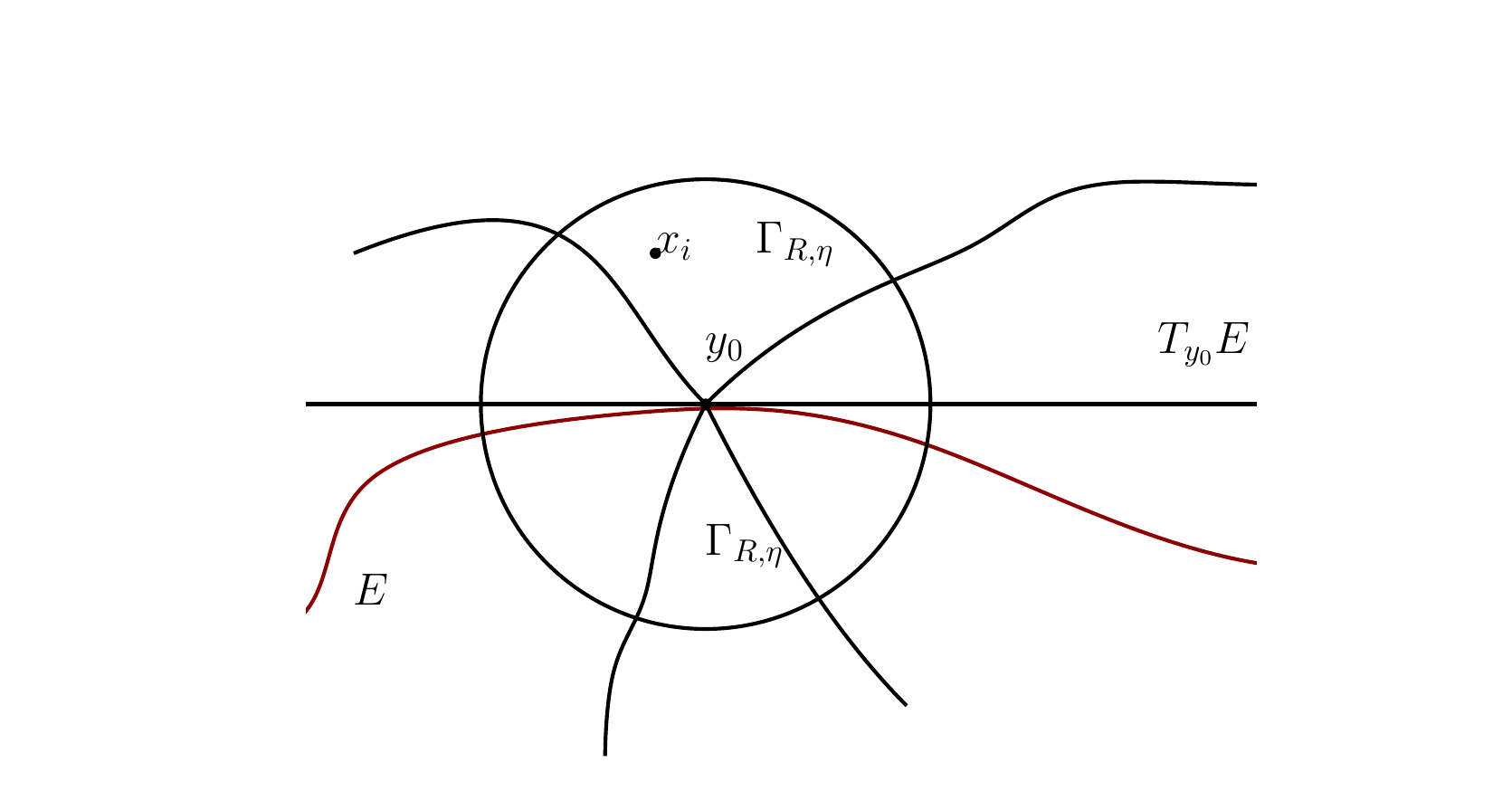}
\end{center}

\begin{thm}\label{NT_theorem}
Let $\mu$ be a $d$-Ahlfors regular measure on $E$ with density $f$ with respect to $\sigma$, and the functions $u_f$ and $\delta$ be as in Section 2. Then, regardless of $\eta$ in the definition of $\Gamma_{R, \eta}$ above, for $\sigma$-almost every $y_0 \in E$ 
\begin{equation}\label{asymp}
    \NTlim_{x \to y_0} u_f(x)\delta(x)^{n - d - 2} = C(n, d)f(y_0),
\end{equation}
where $C(n, d)$ is a constant which depends only on $n$ and $d$.
\end{thm} 

\begin{defn}
We call a measure $\nu$ ($d$-)flat if it is equal to a measure $cdy$, where $c$ is a constant, and $dy$ is a Lebesgue measure supported on a hyperplane of dimension $d$. 
\end{defn}

\begin{rem}\label{ae}
More precisely, in order for the non-tangential limit $\eqref{asymp}$ to exist at a point $y_0$ of $E$, it is enough that for this point $y_0$ the following holds simultaneously: 
\begin{enumerate}
    \item $y_0$ is a density point of the measure $\mu$,
    \item $T_{y_0}E$ exists,
    \item every tangent measure of $\mu$ at $y_0$ is a flat measure.
\end{enumerate}
All three conditions hold, as mentioned above in the statement of the theorem, for almost every point in our set $E$: see $\cite{M}$, the precise reference for the property 2 is given above, and for the property 3 -- in the proof of the theorem. 
\end{rem}

\begin{proof}[Proof of theorem 5]
Our main instruments are the so-called blow-up limits. We fix a point $y_0$ as in Remark \ref{ae} above and a decreasing to zero sequence of positive numbers $\{r_i\}$. Then we introduce the sequence of sets $\{E_i\}$, $E_i = \frac{E - y_0}{r_i}$ and the measures $\mu_i$ supported on $E_i$, $\mu_i(S) = \frac{\mu(r_iS + y_0)}{r_i^d}$. Theorems 14.3 and 16.5 in \cite{M} combined give us that the (weak) limit of (a subsequence of) $\{\mu_i\}$ is a flat tangent measure $\mu_\infty = f(y_0)d\Hd^d = f(y_0)dy$, supported on the hyperplane $E_\infty$, a limit of $\{E_i\}$ in the sense of Hausdorff distance, which coincides as an element of the Grassmannian with the tangent hyperplane $T_{y_0}E$. 

Let us also fix an aperture $0 < \eta < 1$ and a point $x \in \Gamma_{E_\infty, \eta}$, where $\Gamma_{E_\infty, \eta} = \{x \in \R^n \setminus E_\infty: \dist(x, E_\infty) \ge \eta |x|\}$. Note that for $i$ large enough $x$ lies well outside of $E_i$. We define now on $\R^n \setminus E_i$ the function
$$R_i(x) = \ilim_{E_i}{\frac{d\mu_i(y)}{|x - y|^{n - 2}}}.$$ 

Our strategy is the following. Denote $x_i = r_ix + y_0 \in \Gamma_{R, \eta}$. First we will explain that the sequence $\{R_i(x)\}_{i \to \infty}$ has a limit proportional to the limit of $u_f(x_i)\delta(x_i)^{n - d - 2}$, as $x_i$ tends to $y_0$, with the ratio $\dist(x, E_\infty)^{(n - d - 2)}$. Roughly, this means that the limit of the sequence $\{R_i(x)\}_{i \to \infty}$ and the non-tangential limit of the function $u_f(\cdot)\delta(\cdot)^{n - d - 2}$ are proportional. Then we will prove that $\{R_i(x)\}_{i \to \infty}$ converges to the product of the constant $C(n, d)f(y_0)$ and the function $\dist(x, E_\infty)^{-(n - d - 2)}$ uniformly in $x$ which are ``far away'' from $E_\infty$. This gives us what we want.

For the first part of the proof, observe that, with the notation $w_i = r_iy + y_0 \in E$ and $x_i = r_ix + y_0 \in \Gamma_{R, \eta}$,
$$\ilim_{E_i}{\frac{d\mu_i(y)}{|x - y|^{n - 2}}} = \ilim_E{\frac{d\mu(w)_i/r_i^d}{\left|\frac{x_i - y_0}{r_i} - \frac{w_i - y_0}{r_i}\right|^{n - 2}}} = \ilim_E{\frac{d\mu(w_i)}{|x_i - w_i|^{n - 2}}}r_i^{n - 2 - d} = u_f(x_i) r_i^{n - d - 2}.$$
Since $x_i$ lies in the non-tangential access region $\Gamma_{R, \eta}$ (of the point $y_0$), the difference between the distance $\delta(x_i)$ of the point $x_i$ to the set $E$ and the product $|x_i - y_0|\sin{(x \wedge E_\infty)}$, where $x \wedge E_\infty$ is the angle between the hyperplane $E_\infty$ and the vector $x$, tends to zero uniformly in $x$ in $\Gamma_{E_\infty, \eta}$ as $i$ tends to infinity. This implies that, given $|x_i - y_0| = r_i|x| = r_i\frac{\dist(x, E_\infty)}{\sin{(x \wedge E_\infty)}}$,
$$R_i(x) = u_f(x_i) r_i^{n - d - 2} = u_f(x_i) \left(\frac{|x_i - y_0|}{|x|}\right)^{n - d - 2}$$ 
$$ = u_f(x_i) \left(\frac{|x_i - y_0|\sin{(x \wedge E_\infty)}}{\dist(x, E_\infty)}\right)^{n - d - 2} = u_f(x_i)\left(\frac{\delta(x_i)}{\dist(x, E_\infty)}\right)^{n - d - 2} + f(x_i),$$
where $f(x_i)$ tends to zero uniformly in $x$ as $i$ tends to infinity.

For the second part, as announced, we first prove that $R_i(x) \to R_\infty(x)$ uniformly in $x$ which lie ``far away'' from $E_\infty$, where 
$$R_\infty(x) = \ilim_{E_\infty}{\frac{d\mu_\infty(y)}{|x - y|^{n - 2}}}.$$

\begin{lm}\label{unicomp}
Let $K$ be a compact set inside $\R^n \setminus E_\infty$. Then functions $R_i$ converge to $R_\infty$ uniformly on $K$.
\end{lm} 
\begin{proof}
We begin with the observation that if $x \in K$, then for $i$ large enough $x$ lies well outside of the set $E_i$. We consider from now on only such indices $i$. The functions $R_i(x)$ are bounded uniformly in $x \in K$ and $i$. Indeed, the measures $\mu_i$ are $d$-Ahlfors regular with uniform in $i$ constants. If we set $y_0 \in E_i: \dist(x, E_i) = \dist(y_0, x)$ and $r = \min_{i \ge i_0}{\dist(K, E_i)}$ (which is separated from zero), then
$$\ilim_{E_i}{\frac{d\mu_i(y)}{|x - y|^{n - 2}}} = \ilim_{B(y_0, r)}{\frac{d\mu_i(y)}{|x - y|^{n - 2}}} + \slim_{k = 1}^\infty{\ilim_{2^kr \le |y - y_0| \le 2^{k + 1}r}{\frac{d\mu_i(y)}{|x - y|^{n - 2}}}}$$
$$ \le \dist(K, E_i)^{-n + 2}r^d + \slim_{k = 1}^\infty{(2^k r)^{-n + 2} C(2^{k + 1}r)^d} \le C r^{-n + 2 + d}.$$
The uniform estimate implies that for every $\varepsilon > 0$ there is $R > 0$ such that 
$$\left|\ilim_{B(0, R)}{\frac{d\mu_i(y)}{|x - y|^{n - 2}}} - \ilim_{E_i}{\frac{d\mu_i(y)}{|x - y|^{n - 2}}}\right| < \varepsilon$$
for every $i = i_0, \dots, \infty$. There exists a smooth function $\phi$ which approximates $\chi_{B(0, R)}$ and is supported inside $B(0, R)$ such that 
$$\left|\ilim_{B(0, R)}{\frac{d\mu_i(y)}{|x - y|^{n - 2}}} - \ilim_{E_i}{\frac{\phi(y)d\mu_i(y)}{|x - y|^{n - 2}}}\right| < \varepsilon.$$
Therefore we have, for the same set of indices,
$$\ilim_{E_i}{\left|\frac{(1 - \phi(y))d\mu_i(y)}{|x - y|^{n - 2}}\right|} < 2\varepsilon.$$
We need this to work with functions of the type $\frac{\phi(y)}{|x - y|^{n - 2}}$ instead of $\frac{1}{|x - y|^{n - 2}}$.

Observe now that functions of the family $\{\frac{\phi(y)}{|x - y|^{n - 2}}\}_{x \in K}$ are still uniformly bounded and moreover equicontinuous inside $\overline{B(0, R)}$. Indeed, if $|x_1 - x_2| < \delta$, then
$$\left|\frac{\phi(y)}{|x_1 - y|^{(n - 2)}} - \frac{\phi(y)}{|x_2 - y|^{n - 2}}\right| \le \frac{\left||x_2 - y|^{n - 2} - |x_1 - y|^{n - 2}\right|}{|x_1 - y|^{n - 2}|x_2 - y|^{n - 2}} $$
$$ \le (|x_2 - y| - |x_1 - y|)\slim_{i = 0}^{n - 3}{\frac{|x_2 - y|^{n - 3 - i}|x_1 - y|^i}{|x_1 - y|^{n - 2}|x_2 - y|^{n - 2}}} \le |x_2 - x_1|C(n) < \delta C(n).$$
Therefore we can apply Arzela-Ascoli theorem to the family $\{\frac{\phi(y)}{|x - y|^{n - 2}}\}_{x \in K}$ and find a finite collection of continuous function $\{g_k\}$ with support inside $B(0, 2R)$ such that for every point $x \in K$ there is an index $k$ with $|g_k(y) - \frac{\phi(y)}{|x - y|^{n - 2}}| \le \varepsilon R^{-d}$ for $y \in B(0, 2R) \cap E_i$. Then
$$\ilim{\left|g_k(y) - \frac{\phi(y)}{|x - y|^{n - 2}}\right|d\mu_i(y)} + \ilim{\left|g_k(y) - \frac{\phi(y)}{|x - y|^{n - 2}}\right|d\mu_\infty(y)} \le C \varepsilon.$$
By the definition of a tangent measure $\ilim{g_k(y)d\mu_i(y)}$ converges to $\ilim{g_k(y)d\mu_\infty(y)}$, and we can finally estimate the modulus of the difference $|R_i(x) - R_\infty(x)|$ the following way:
$$\left|\ilim_{E_i}{\frac{d\mu_i(y)}{|x - y|^{n - 2}}} - \ilim_{E_\infty}{\frac{d\mu_\infty(y)}{|x - y|^{n - 2}}}\right| \le 4\varepsilon + \left|\ilim_{E_i}{\frac{\phi(y)d\mu_i(y)}{|x - y|^{n - 2}}} - \ilim_{E_\infty}{\frac{\phi(y)d\mu_\infty(y)}{|x - y|^{n - 2}}}\right|$$
$$\le 4\varepsilon + \ilim_{E_i}{\left|g_k(y) - \frac{\phi(y)}{|x - y|^{n - 2}}\right|d\mu_i(y)} + \ilim_{E_\infty}{\left|g_k(y) - \frac{\phi(y)}{|x - y|^{n - 2}}\right|d\mu_\infty(y)}$$ 
$$+ \left|\ilim{g_k(y)d\mu_i(y)} - \ilim{g_k(y)d\mu_\infty(y)}\right| \le C\varepsilon.$$
\end{proof}

We now know that, on the one hand, the limit of the function $R_i(x)$, as $i$ tends to infinity, coincides with the limit of the function
$$u_f(x_i)\delta(x_i)^{n - d - 2}\dist(x, E_\infty)^{-(n - d - 2)}$$
as $x_i$ tends to $y_0$, if $x_i = r_i x + y_0$. On the other hand, $R_i(x)$ tends to $R_\infty(x)$ uniformly for $x$ in $\Gamma_{E_\infty, \eta}$ which stay away from $E_\infty$. This implies that for every $x$ such that $|x| = 1$, say, and $x \in \Gamma_{E_\infty, \eta}$ 
$$u_f(x_i)\delta(x_i)^{n - d - 2} \to R_\infty(x)\dist(x, E_\infty)^{n - d - 2}, \; i \to \infty, x_i = r_i x + y_0.$$
The product $R_\infty(x)\dist(x, E_\infty)^{n - d - 2}$ we will compute below, and it is a product of $C(n, d)f(y_0)$, where the constant $C(n, d)$ can be calculated explicitly in terms of $\Gamma$-functions and depends only on $d$ and $n$. This is exactly the right-hand side of $\eqref{asymp}$. It is left for us only to see why
$$u_f(x_i)\delta(x_i)^{n - d - 2} \to R_\infty(x)\dist(x, E_\infty)^{n - d - 2}, \quad \mbox{as} \; x_i \to y_0$$
for any sequence $\{x_i\}$ in $\Gamma_{R, \eta}$, not just those of the type $x_i = r_ix + y_0$ for a fixed $x$. But this follows from the uniform convergence of 
$$R_i(x) \quad \mbox{to} \quad u_f(x_i)\left(\frac{\delta(x_i)}{\dist(x, E_\infty)}\right)^{n - d - 2}$$
and in Lemma \ref{unicomp}: we have $|u_f(x_i)\delta(x_i)^{n - d - 2} - C(n, d)f(y_0)| < \epsilon$ as soon as $|x_i - y_0|$ is small enough. 

~\

To finish the proof we compute $R_\infty(x)\dist(x, E_\infty)^{n - d - 2}.$ By definition,
$$R_\infty(x)\dist(x, E_\infty)^{n - d - 2} = \ilim_{E_\infty}{\frac{f(y_0)dy}{|x - y|^{n - 2}}}\dist(x, E_\infty)^{n - d - 2}.$$
Let $z$ be the point of $E_\infty$ such that $\dist(x, E_\infty) = \dist(x, z)$. Clearly $|x - y| = \sqrt{|x - z|^2 + r^2},$ where $r = \dist(z, y)$ for $y \in E_\infty$. Write $\delta_1(x)$ for $\dist(x, z)$. Then we have
$$\ilim_{E_\infty}{\frac{dy}{|x - y|^{n - 2}}} = \ilim_0^\infty{\frac{r^{d - 1}dr}{(\delta_1(x)^2 + r^2)^{(n - 2)/2}}}$$ 
$$ = \delta_1(x)^{-n + d + 2} \ilim_0^\infty{\frac{\frac{r^{d - 1}}{\delta_1(x)^{d - 1}} d\left(\frac{r}{\delta_1(x)}\right)}{\left(1 + \left(\frac{r}{\delta_1(x)}\right)^2\right)^{(n - 2)/2}}} = \delta_1(x)^{-n + d + 2}\ilim_0^\infty{\frac{x^{d - 1}dx}{(1 + x^2)^{(n - 2)/2}}}.$$

Note that the integral $\ilim_0^\infty{\frac{x^{d - 1}dx}{(1 + x^2)^{(n - 2)/2}}}$ clearly converges and is equal to $c_1 = C(n, d) = \Vol(\Sf^{d - 1})\frac{1}{2}\frac{\Gamma\left(\frac{d}{2}\right)\Gamma\left(\frac{n - d - 2}{2}\right)}{\Gamma\left(\frac{n - 2}{2}\right)}$. Therefore $R_\infty(x)\dist(x, E_\infty)^{n - d - 2}$ is indeed equal to $C(n, d)f(y_0)$.

\end{proof}

\begin{cor}\label{cor1}
Apart from the fact that the function 
$$u_f(x)\delta(x)^{n - d - 2} = \delta(x)^{n - d - 2}\ilim_{E}{\frac{f(y)dy}{|x - y|^{n - 2}}},$$
where $f \in L^\infty(E)$, has a non-tangential limit at $\sigma$-almost every $y_0 \in E$, which is equal to $c_1f(y_0)$, throughout the next sections we will also use similar statements, but for degrees in the denominator in the integral other than $(n - 2)$. More precisely, we claim that for $f \in L^\infty(E)$ and $\beta > 0$ the function 
$$\delta(x)^{\beta}\ilim_{E}{\frac{f(y)dy}{|x - y|^{d + \beta}}},$$ 
also has a non-tangential limit at almost every $y_0 \in E$ equal to $C(d, \beta)f(y_0)$, and the constant $C(d, \beta)$ depends only on $d$ and $\beta$. The proof is similar to the proof of $\eqref{asymp}$.
\end{cor}

\section{Three BMO lemmas}

In this section we prove some more technical preliminaries concerning the space of functions of bounded mean oscillations (BMO) to simplify the exposition of Sections 5 and 6. The reader can skip to them directly and return to this section when needed.

By $x$ we usually (that is, except for in Lemma $\ref{lmaverages}$) denote a point in $\R^n$ which lies away from the hyperplane $\R^d$ containing zero. By $\delta(x)$ we denote the distance between $x$ and the hyperplane $\R^d$. 

Recall that the John-Nirenberg inequality asserts that for every $1 \le p < \infty$
$$\sup_B{\frac{1}{|B|}\left(\ilim_B{|f(y) - m_Bf|^p dy}\right)^{1/p}} \asymp \|f\|_{BMO},$$
where by $m_Bf$ we denote the average $\frac{1}{|B|}\ilim_B{f(y)dy}$, and the supremum is taken over all balls in $\R^d$. 

\begin{lm}\label{lmaverages}
Let $f$ be a function in $BMO(\R^d)$. Denote by $B(x, r)$ the ball with center $x$ and radius $r$. Then for every $r > r_0 > 0$ holds
\begin{equation}\label{averages}
    |m_{B(x, r)}f - m_{B(x_0, r_0)}f| \le C\|f\|_{BMO}\left(\ln{\frac{r}{r_0}} + \ln{\frac{|x - x_0|}{r_0}}\right).   
\end{equation}
\end{lm}
\begin{proof}
We will need two simple facts. First, let $B$ and $B'$ be balls such that $B' \subset B$. Then
\begin{equation}\label{av1}
    |m_Bf - m_{B'}f| \le \frac{|B|}{|B'|}\|f\|_{BMO}.
\end{equation}
Indeed,
$$|m_Bf - m_{B'}f| = \Big|\frac{1}{|B'|}\ilim_{B'}{f(y)dy} - m_{B}f\Big| \le \frac{1}{|B'|}\ilim_{B'}{|f(y) - m_{B}f|}$$ $$ \le \frac{|B|}{|B'|}\frac{1}{|B|}\ilim_B{|f(y) - m_{B}f|} \le \frac{|B|}{|B'|}\|f\|_{BMO}.$$
Second, if the two centers $x$ and $x'$ are such that $|x - x'| \le 2r$, then
\begin{equation}\label{av2}
    |m_{B(x, r)}f - m_{B(x', r)}f| \le C\|f\|_{BMO}.
\end{equation}
This follows from $\eqref{av1}$: we can find a ball $B$ of radius $3r$ such that $B(x, r) \subset B$ and $B(x', r) \subset B$. Then
$$|m_{B(x, r)}f - m_{B(x', r)}f| \le |m_{B(x, r)}f - m_Bf| + |m_Bf - m_{B(x', r)}f| \le C 3^d\|f\|_{BMO}.$$

We are now ready to prove $\eqref{averages}$. 
$$|m_{B(x, r)}f - m_{B(x_0, r_0)}f| \le |m_{B(x, r)}f - m_{B(x, r_0)}f| + |m_{B(x, r_0)}f - m_{B(x_0, r_0)}f|.$$
Inequality $\eqref{av1}$ will help us to estimate the first term. Build a chain of nested balls $B_1 = B(x, r_0) \subset B_2 \subset \dots \subseteq B_N = B(x, r)$ with center $x$ and radii $r_i$ such that $r_i = 2r_{i - 1}$ for $i \le N - 1$. Then it is clear that $N \le C\ln{\frac{r}{r_0}}$. For each $1 < i \le N$ we have $|m_{B_{i - 1}}f - m_{B_i}f| \le C\|f\|_{BMO}$. Therefore
$$|m_{B(x, r)}f - m_{B(x, r_0)}f| \le \slim_{i = 2}^N{|m_{B_i}f - m_{B_{i - 1}}f|} \le C\|f\|_{BMO}\ln{\frac{r}{r_0}}.$$
Inequality $\eqref{av2}$ will help us to estimate the second term $|m_{B(x, r_0)}f - m_{B(x_0, r_0)}f|$. Build a chain of balls $B_1 = B(x, r_0), B_2, \dots, B_N = B(x_0, r_0)$ of common radii $r_0$ and centers $x_i$ such that for $i \le N - 1$ we have $|x_{i - 1} - x_i|  = 2r_0$. Clearly $N \le C \ln{\frac{|x - x_0|}{r_0}}$. By $\eqref{av2}$,
$$|m_{B(x, r_0)}f - m_{B(x_0, r_0)}f| \le \slim_{i = 2}^N{|m_{B_i}f - m_{B_{i - 1}}f|} \le C\|f\|_{BMO}\ln{\frac{|x - x_0|}{r_0}}.$$
This completes the proof of $\eqref{averages}$.
\end{proof}

\begin{lm}
Let $f$ be a function in $BMO (\R^d)$ and $x$ be a point in $\R^n \setminus \R^d$. Then for any integer $m \ge 1$
\begin{equation}\label{f_bmo3} 
    \ilim_{\R^d}{\frac{|f(y)|^m dy}{|x - y|^{d + \beta}}} \le C(1 + \delta(x)^{-(d + \beta)}) \left(\|f\|_{BMO}^m + |m_{B(x_0, 1)}f|^m\right),
\end{equation}
where by $x_0$ we denote the projection of $x$ to the hyperplane $\R^d$ ($x_0 = (x_1, \dots , x_d)$), and the constant $C$ depends only on $d, m$ and $\beta > 0$.
\end{lm}
\begin{proof}
First, observe that
$$\ilim_{\R^d}{\frac{|f(y)|^m dy}{|x - y|^{d + \beta}}} \le C(m) \ilim_{\R^d}{\frac{|f(y) - m_{B(x_0, 1)}f|^m dy}{|x - y|^{d + \beta}}} + C(m) \ilim_{\R^d}{\frac{|m_{B(x_0, 1)}f|^m dy}{|x - y|^{d + \beta}}}.$$
We use a computation very similar to the one we saw at the end of the previous section to calculate $\ilim_{\R^d}{|x - y|^{-(d + \beta)}dy}$ and get
\begin{equation}\label{lm6_1}
    \ilim_{\R^d}{\frac{|m_{B(x_0, 1)}f|^m dy}{|x - y|^{d + \beta}}} \le C\delta(x)^{- \beta}|m_{B(x_0, 1)}f|^m \le C(1 + \delta(x)^{-(d + \beta)})|m_{B(x_0, 1)}f|^m.
\end{equation}
The second inequality here is true because $\delta(x)^{-\beta} \le \delta(x)^{-(d + \beta)}$ if $\delta(x) \le 1$, and otherwise both quantites are dominated by a constant. This means that we can assume $m_{B(x_0, 1)f} = 0$. 

Without loss of generality, we can also assume that the fist $d$ and the last $n - d - 1$ coordinates of $x$ are zero: $x = (0, \dots, 0, \delta(x), 0, \dots)$. We split the space $\R^d$ into $B(1, 0)$ and the union of rings $B(0, 2^k) \setminus B(0, 2^{k - 1})$ and estimate separately integrals of $|f(y)|^m/|x - y|^{d + \beta}$ over those sets. For the first term we have
$$\ilim_{B(0, 1)}{\frac{|f(y)|^m dy}{|x - y|^{d + \beta}}} = \ilim_{B(0, 1)}{\frac{|f(y) - m_{B(0, 1)f}|^m dy}{|x - y|^{d + \beta}}} $$ 
\begin{equation}\label{lm6_2}
    \le C\delta(x)^{- (d + \beta)}\frac{1}{|B(0, 1)|}\ilim_{B(0, 1)}{|f(y) - m_{B(0, 1)f}|^mdy} \le C\delta(x)^{- (d + \beta)}\|f\|_{BMO}^m
\end{equation}
by the John-Nirenberg inequality. Then we estimate the integrals over the rings; for $B = B(0, 2^{k - 1})$ and $2B = B(0, 2^k)$ we have
$$\ilim_{2B \setminus B}{\frac{|f(y)|^m dy}{|x - y|^{d + \beta}}} \le C(m)\ilim_{2B \setminus B}{\frac{|f(y) - m_{2B}f|^m dy}{|x - y|^{d + \beta}}} + C(m)\ilim_{2B \setminus B}{\frac{|m_{2B}f|^m dy}{|x - y|^{d + \beta}}}.$$
Since for $y$ in $2B \setminus B$ the quantity $|x - y|^2$ is comparable to $\delta(x)^2 + 2^{2k}$, $|x - y|^{-(d + \beta)}$ is dominated by $C(\delta(x) + 2^k)^{-(d + \beta)}$. For the first term this gives
$$\ilim_{2B \setminus B}{\frac{|f(y) - m_{2B}f|^m dy}{|x - y|^{d + \beta}}} \le C(\delta(x) + 2^k)^{- (d + \beta)}2^{kd}\fint_{2B}{|f(y) - m_{2B}f|^m dy} $$ 
$$ \le C(\delta(x) + 2^k)^{- (d + \beta)}2^{kd}\|f\|_{BMO}^m,$$
again by John-Nirenberg. For the second term we have
$$\ilim_{2B \setminus B}{\frac{|m_{2B}f|^m dy}{|x - y|^{d + \beta}}} = \ilim_{2B \setminus B}{\frac{|m_{2B}f - m_{B(0, 1)}f|^m dy}{|x - y|^{d + \beta}}} \le C 2^{kd} (\delta(x) + 2^k)^{- (d + \beta)}|m_{2B}f - m_{B(0, 1)}f|^m$$ 
$$\le 2^{kd} (\delta(x) + 2^k)^{- (d + \beta)} (\ln{2^k})^m \|f\|_{BMO}^m \le C2^{kd} (\delta(x) + 2^k)^{- (d + \beta)} k^m \|f\|_{BMO}^m$$
by $\eqref{averages}$. Therefore
$$\ilim_{2B \setminus B}{\frac{|f(y)|^m dy}{|x - y|^{d + \beta}}} \le C 2^{kd} (\delta(x) + 2^k)^{- (d + \beta)} (1 + k^m)\|f\|_{BMO}^m.$$
Obsesrve that the sum $\slim_{k = 1}^\infty {2^{kd} (\delta(x) + 2^k)^{- (d + \beta)} (1 + k^m)}$ is finite and dominated by a universal constant $C$, since $(\delta(x) + 2^k)^{-1} \le 2^{-k}$. Then we sum over all the rings and get
\begin{equation}\label{lm6_3}
    \ilim_{\R^d \setminus B(0, 1)}{\frac{|f(y)|^m dy}{|x - y|^{d + \beta}}} = \slim_{k = 1}^\infty {\ilim_{B(0, 2^k) \setminus B(0, 2^{k - 1})}{\frac{|f(y)|^m dy}{|x - y|^{d + \beta}}}} \le C\|f\|_{BMO}^m.
\end{equation}
Collecting $\eqref{lm6_1}$, $\eqref{lm6_2}$ and $\eqref{lm6_3}$, we obtain the final estimate $\eqref{f_bmo3}$.
\end{proof}

\begin{cor}
For the same notation as in $\eqref{f_bmo3}$, the estimate
\begin{equation}\label{f_bmo4}
    \ilim_{\R^d}{\frac{|f(y)|^m dy}{|x - y|^{d + \beta}}} \le C (1 + \delta(x)^{-(d + \beta)})(1 + \ln{|x_0|})^m\left(\|f\|_{BMO}^m + |m_{B(0, 1)}f|^m\right)
\end{equation}
holds.
\end{cor}
\begin{proof}
This follows from $\eqref{f_bmo3}$ and Lemma $\ref{lmaverages}$.
\end{proof}

\begin{lm}
Let $f_i, i = 1, \dots, m$ be functions in $BMO(\R^d)$. With the same notation as in $\eqref{f_bmo3}$, we have
\begin{equation}\label{Ho}
    \ilim_{\R^d}{\frac{|f_1(y) \dots f_m(y)|dy}{|x - y|^{d + \beta}}} \le C( 1 + \delta(x)^{-(d + \beta)})(1 + \ln{|x_0|})^m \prod_{i = 1}^m{(\|f_i\|_{BMO}^m + |m_{B(0, 1)}f_i|^m)^{1/m}}
\end{equation}
\end{lm}
\begin{proof}
This is, essentially, the H\"older inequality applied consequentially $m$ times. We give a proof by induction for the sake of completeness. First, we want to prove that
\begin{equation}\label{aux_Ho}
    \ilim_{\R^d}{\frac{|f_1(y) \dots f_m(y)|dy}{|x - y|^{d + \beta}}} \le \prod_{i = 1}^m{\left(\ilim_{\R^d}{\frac{|f_i(y)|^m dy}{|x - y|^{d + \beta}}}\right)^{1/m}}.
\end{equation}
For the base we apply H\"older with exponents $p_1 = m$ and $q_1 = \frac{m - 1}{m}$ and functions $|f_1(y)||x - y|^{-(d + \beta)\frac{1}{m}}$ and $|f_2(y) \dots f_m(y)||x - y|^{-(d + \beta)\frac{m - 1}{m}}$. It gives
$$\ilim_{\R^d}{\frac{|f_1(y) \dots f_m(y)|dy}{|x - y|^{d + \beta}}} \le \left(\ilim_{\R^d}{\frac{|f_1(y)|^m dy}{|x - y|^{d + \beta}}}\right)^{1/m}\left(\ilim_{\R^d}{\frac{|f_2(y) \dots f_m(y)|^\frac{m}{m - 1} dy}{|x - y|^{d + \beta}}}\right)^{(m - 1)/m}.$$
By the induction hypothesis,
$$\ilim_{\R^d}{\frac{|f_1(y) \dots f_m(y)|dy}{|x - y|^{d + \beta}}} \le \prod_{i = 1}^j{\left(\ilim_{\R^d}{\frac{|f_i(y)|^m dy}{|x - y|^{d + \beta}}}\right)^{1/m}}\left(\ilim_{\R^d}{\frac{|f_{j + 1}(y) \dots f_m(y)|^\frac{m}{m - j} dy}{|x - y|^{d + \beta}}}\right)^{(m - j)/m}.$$
If $j = m - 1$, we are done, otherwise we apply H\"older once again to 
$$\ilim_{\R^d}{\frac{|f_{j + 1}(y) \dots f_m(y)|^\frac{m}{m - j} dy}{|x - y|^{d + \beta}}}$$
with $p_{j + 1} = m - j$, $q_{j + 1} = \frac{m - (j + 1)}{m - j}$ and functions 
$$|f_{j + 1}(y)|^\frac{m}{m - j}|x - y|^{-(d + \beta)\frac{1}{p_{j + 1}}}, \quad |f_{j + 2}(y) \dots f_m(y)|^\frac{m}{m - j}|x - y|^{-(d + \beta)\frac{1}{q_{j + 1}}}.$$ 
We will get
$$\ilim_{\R^d}{\frac{|f_1(y) \dots f_m(y)|dy}{|x - y|^{d + \beta}}} \le \prod_{i = 1}^{j + 1}{\left(\ilim_{\R^d}{\frac{|f_i(y)|^m dy}{|x - y|^{d + \beta}}}\right)^{1/m}}\left(\ilim_{\R^d}{\frac{|f_{j + 2}(y) \dots f_m(y)|^\frac{m}{m - (j + 1)} dy}{|x - y|^{d + \beta}}}\right)^{(m - (j + 1))/m},$$
which completes the induction step.

The inequality $\eqref{Ho}$ follows from $\eqref{f_bmo4}$ and $\eqref{aux_Ho}$.
\end{proof}

\section{No one-parameter families of solutions}

We can now turn to the study of solutions of the equation $\eqref{solution}$ for rectifiable sets $E$ with integer dimension $d$, using our preliminaries from sections 2 and 3. Assume the function $D_{\alpha, \mu}$ as in $\eqref{mu_dist}$ is a solution of $\eqref{solution}$. Then we know that, see the beginning of section 2, for $\gamma = - n + d + 2$, $D_{\alpha, \mu}^\gamma$ is equivalent to $\delta(x)^{-n + d + 2}$ and is harmonic outside $E$. Therefore, according to subsection 2.2, there exists a density function $h$ in $L^\infty(E)$ such that
\begin{equation}\label{functions}
    \ilim_{E}{\frac{h(y)d\sigma(y)}{|x - y|^{n - 2}}} =  \left(\ilim_{E}{\frac{f(y)d\sigma(y)}{|x - y|^{d + \alpha}}}\right)^{(n - d - 2)/\alpha} \quad \forall \; x \in \R^n \setminus E,
\end{equation}
where $f \in L^\infty(E)$ is just the density of $\mu$ with respect to $\sigma = \Hd^d|_E$.
Multiplying both sides of $\eqref{functions}$ by $\delta(x)^{n - d - 2}$ and using $\eqref{asymp}$ in Theorem 5 and Corollary \ref{cor1} to pass to non-tangential limits, we get that for almost every $y \in E$
\begin{equation}\label{dens_relation}
    c_1h(y) = (c_2 f(y))^{(n - d - 2)/\alpha},
\end{equation}
where the constants $c_i$ depend only on $n, d$ and $\alpha$. The constant $c_1$ is the constant $C(n, d)$ in $\eqref{asymp}$, which we computed at the end of the proof of theorem 5, and the constant $c_2$ is equal to $\delta(x)^{-\alpha}\ilim_{\R^d}{}\frac{dy}{|x - y|^{d + \alpha}} = \Vol(\Sf^{d - 1})\frac{1}{2}\frac{\Gamma\left(\frac{d}{2}\right)\Gamma\left(\frac{\alpha}{2}\right)}{\Gamma\left(\frac{d + \alpha}{2}\right)}$. Thus, we get an equation for the density $f$ of the measure $\mu$, when $D_{\alpha, \mu}$ satisfies $\eqref{solution}$, but we prefer to write everything in terms of the density function $h$ alone. So, combining $\eqref{functions}$ and $\eqref{dens_relation}$, we get that $D_{\alpha, \mu}$ being a solution of $\eqref{solution}$ is equivalent to the equation
\begin{equation}\label{densities}
    \ilim_{E}{\frac{h(y)dy}{|x - y|^{n - 2}}} = c_3 \left(\ilim_{E}{\frac{h(y)^\frac{\alpha}{n - d - 2}dy}{|x - y|^{d + \alpha}}}\right)^\frac{n - d - 2}{\alpha}, \quad \forall \; x \in \R^n \setminus E,
\end{equation}
where we denote by $c_3$ the constant $c_2^{-\frac{n - d - 2}{\alpha}} c_1$. Note that the density function $h$ is also bounded away from zero, since according to Lemma $\ref{mes2}$ the density function $f$ has this property, and the two functions are connected by $\eqref{dens_relation}$.

~\

It is easy to check that for $E = \R^d$ or any other hyperplane a constant function $h$ gives a solution for $\eqref{densities}$. This is what we call the flat solution (because the measure $c\Hd^d$ on $\R^d$ is flat). As we have said in the introduction, we would be happy to show that no other solutions exist in a small neighbourhood of the flat solution. This is the same as saying that there is no family of solutions with elements arbitrarily close to the flat solution. Our strategy is to linearise $\eqref{densities}$, or to take a derivative in some sense at constant function $h$ and $E = \R^d$: for an arbitrary family of solutions it would be the Frechet derivative. It turns out though, it does not work for arbitrary family of solutions with no regularity. Even when we restrict ourselves to the flat case $E = \R^d$ and are trying to figure out if there is a neighbourhood of the flat solution for which no densities $h$ solve $\eqref{densities}$. The reason for it is, morally, the following. If we suppose the contrary, we can find arbitrarily close to $h = 1$ (or to any other constant) a density $1 + \phi_t$ for which
$$\ilim_{\R^d}{\frac{(1 + \phi_t)dy}{|x - y|^{n - 2}}} = c_3 \left(\ilim_{\R^d}{\frac{(1 + \phi_t)^\frac{\alpha}{n - d - 2}dy}{|x - y|^{d + \alpha}}}\right)^\frac{n - d - 2}{\alpha}, \quad \forall \; x \in \R^n \setminus \R^d,$$
where the parameter $t$ denotes, vaguely, the size of the neighbourhood to which $\phi_t$ belongs. Then we could use the Banach-Alaoglu theorem in a suitable functional space as a compactness argument to say that the family of functions $\{\phi_t\}$ normalized correctly has a non-constant weak limit. This weak limit would satisfy a linear convolution equation, corresponding to the fact that the Fr\'echet derivative of $\eqref{densities}$ is equal to zero ``at the constant solution''. This equation, as we will see later, has no solutions except for constants. Which would finish the proof by contradiction. However, the right space for the functions $\phi_t$ for the compactness argument to work turns out to be the homogeneous Besov space $\dot B^0_{\infty, \infty}$. For this space it seems that the Fr\'echet derivative does not exist. That is, we have a natural space for the compactness argument ($\dot B^0_{\infty, \infty}$), and a smaller one (for instance, BMO) for the existence of a derivative. See subsection $5.3$ for the further comments about this.

Fortunately, our strategy still works if we introduce a parameter $t$ and add the assumption of differentiability at zero of the family $\{\phi_t\}$ in $BMO$ with respect to it and a reasonable boundedness assumption. This is exactly what the assumptions 2 and 4 in Definitions $\ref{defn1}$ and $\ref{defn2}$ are about. We are going to implement our strategy with the indicated assumptions in the next subsections.

\subsection{No one-parameter families of flat solutions}

In this subsection we prove Theorem $\ref{mainthm_flat}$. Suppose that there is a non-trivial one-parameter family of perturbations in the plane $\R^d$ of the flat solution of $\eqref{densities}$ with the density $h = 1$. See Definition \ref{defn1} to recall what it is. Note that, even though the densities $1 + \phi_t$ are in $L^\infty(\R^d)$ and therefore $\phi_t$ is in the space $L^\infty \cap BMO (\R^d) = L^\infty(\R^d)$, it is still more natural to consider $\phi_t$ as an element of $BMO (\R^d)$ (or as an equivalence class in $L^\infty(\R^d)/\sim$, where $f \sim g$ if $f - g$ is a constant function), because $1 + \phi_t(y) + c$ is a density function of a perturbation of the flat solution with the density $h_c = 1 + c$. We can assume that $\|\phi_t\|_\infty < 1/2$. 

For every $t$ the equation $\eqref{densities}$ for the solution from our family $D_{\alpha, \mu_t}$ with the density $1 + \phi_t(y)$ turns into
\begin{equation}\label{take_diff}
    I(t, x) := \ilim_{\R^d}{\frac{(1 + \phi_t)dy}{|x - y|^{n - 2}}} - c_3\left(\ilim_{\R^d}{\frac{(1 + \phi_t)^\frac{\alpha}{n - d - 2}dy}{|x - y|^{d + \alpha}}}\right)^\frac{n - d - 2}{\alpha} = 0, \quad \forall \; x \in \R^n \setminus \R^d.
\end{equation}
We want to differentiate the left-hand side of $\eqref{take_diff}$ at $t = 0$. The natural candidate for the derivative at zero would be, of course,
$$
    \frac{\partial I(t, x)}{\partial t}(0) = \ilim_{\R^d}{\frac{\frac{\phi_t}{\partial t}(0, y)dy}{|x - y|^{n - 2}}} - c_3\left(\ilim_{\R^d}{\frac{(1 + \phi_{0})^\frac{\alpha}{n - d - 2}dy}{|x - y|^{d + \alpha}}}\right)^{\frac{n - d - 2}{\alpha} - 1}\ilim_{\R^d}{\frac{(1 + \phi_{0})^{\frac{\alpha}{n - d - 2} - 1}\frac{\partial \phi_t}{\partial t}(0, y)dy}{|x - y|^{d + \alpha}}},
$$
so, by the definition of $c_2$, and since $\phi_0 = 0$,
\begin{equation}\label{t'}
    \frac{\partial I(t, x)}{\partial t}(0) = \ilim_{\R^d}{\frac{\frac{\phi_t}{\partial t}(y, 0)dy}{|x - y|^{n - 2}}} - c_1c_2^{-1}\delta(x)^{d + 2 + \alpha - n}\ilim_{\R^d}{\frac{\frac{\partial \phi_t}{\partial t}(y, 0)dy}{|x - y|^{d + \alpha}}}.
\end{equation}
Note that here we use $\frac{\partial I(t, x)}{\partial t}(0)$ just as a notation. To show that $\eqref{t'}$ is indeed the Frechet derivative of $I(t, x)$ we prove the following lemma. This is when our assumptions 2 and 4 from Definition 1 of the one-parameter family of perturbations come into play.

\begin{lm}
For every $x$ in $\R^n \setminus \R^d$ 
$$\left|I(t, x) - I(0, x) - t\frac{\partial I(t, x)}{\partial t}(0)\right| =  \left|I(t, x) - t\frac{\partial I(t, x)}{\partial t}(0)\right|  = o(t), \quad t \to 0.$$
\end{lm}
\begin{proof}
We need to prove that
$$\left|I(t, x) - t\frac{\partial I(t, x)}{\partial t}(0)\right| = \bigg| \ilim_{\R^d}{\frac{(1 + \phi_t)dy}{|x - y|^{n - 2}}} - t\ilim_{\R^d}{\frac{\frac{\phi_t}{\partial t}(y, 0)dy}{|x - y|^{n - 2}}} + \bigg.
$$
\begin{equation}\label{F}
\bigg. + c_1c_2^{-1}\delta(x)^{d + 2 + \alpha - n}t\ilim_{\R^d}{\frac{\frac{\partial \phi_t}{\partial t}(y, 0)dy}{|x - y|^{d + \alpha}}} - c_3\left(\ilim_{\R^d}{\frac{(1 + \phi_t)^\frac{\alpha}{n - d - 2}dy}{|x - y|^{d + \alpha}}}\right)^\frac{n - d - 2}{\alpha} \bigg| = o(t), \quad t \to 0.
\end{equation}
Denote 
$$I_1(t, x) = \ilim_{\R^d}{\frac{(1 + \phi_t)dy}{|x - y|^{n - 2}}} - t\ilim_{\R^d}{\frac{\frac{\partial \phi_t}{\partial t}(y, 0)dy}{|x - y|^{n - 2}}}, \quad I_2(t, x) = \left(\ilim_{\R^d}{\frac{(1 + \phi_t)^\frac{\alpha}{n - d - 2}dy}{|x - y|^{d + \alpha}}}\right)^\frac{n - d - 2}{\alpha}, \quad \mbox{and}$$
$$I_3(t, x) = c_1c_2^{-1}\delta(x)^{d + 2 + \alpha - n}t\ilim_{\R^d}{\frac{\frac{\partial \phi_t}{\partial t}(y, 0)dy}{|x - y|^{d + \alpha}}} - c_3\left(\ilim_{\R^d}{\frac{(1 + \phi_t)^\frac{\alpha}{n - d - 2}dy}{|x - y|^{d + \alpha}}}\right)^\frac{n - d - 2}{\alpha}.$$
\begin{cl}
\begin{equation}\label{I_2}
    I_2(t, x) = c_2^\frac{n - d - 2}{\alpha}\delta(x)^{d + 2 - n} + c_2^{\frac{n - d - 2}{\alpha} - 1}\delta(x)^{d + 2 + \alpha - n}\ilim_{\R^d}{\frac{\phi_t(y)dy}{|x - y|^{d + \alpha}}} + o(t), \; t \to 0.
\end{equation}
\end{cl}
\begin{proof}
This is, essentially, a simple computation which uses Taylor's theorem. Recall that $\|\phi_t\|_\infty \le 1/2 \; \forall t$. First, we rewrite the nominator inside the integral in $I_3$ as
$$(1 + \phi_t(y))^\frac{\alpha}{n - d - 2} = 1 + \frac{\alpha}{n - d - 2}\phi_t(y) + \phi_t^2(y)g_t(y),$$
where $g_t(y)$ is a bounded uniformly in $t$ function. Here we just used an observation that for a real number $x$ such that $|x| \le 1/2$ we have $(1 + x)^\gamma = 1 + \gamma x + c x^2$, where $c$ is a constant uniformly bounded in $x$. Second, we rewrite $I_3$ itself as
$$\left(\ilim_{\R^d}{\frac{dy}{|x - y|^{d + \alpha}}} + \ilim_{\R^d}{\frac{\frac{\alpha}{n - d - 2}\phi_t(y) + \phi_t(y)^2g_t(y)}{|x - y|^{d + \alpha}}}\right)^\frac{n - d - 2}{\alpha}$$ $$ = c_2^\frac{n - d - 2}{\alpha}\delta(x)^{d + 2 - n} + c_2^{\frac{n - d - 2}{\alpha} - 1}\delta(x)^{d + 2 + \alpha - n}\ilim_{\R^d}{\frac{\phi_t(y)dy}{|x - y|^{d + \alpha}}} + \ilim_{\R^d}{\frac{\phi_t(y)^2 G_t(y)dy}{|x - y|^{d + \alpha}}},$$
where $G_t(y)$ is also a bounded uniformly in $t$ function. It is left to show that the last term is $o(t)$. But $\eqref{f_bmo4}$ applied with $m = 2$, combined with boundedness of $G_t(y)$, gives that
$$\ilim_{\R^d}{\frac{\phi_t(y)^2 G_t(y)dy}{|x - y|^{d + \alpha}}} \le C(1 + \delta(x)^{- (d + \alpha)})(\|\phi_t\|_{BMO}^2 + |m_{B(0, 1)}\phi_t|^2).$$
Condition 4 in Definition $\ref{defn1}$ guarantees that $\|\phi_t\|_{BMO}^2 + |m_{B(0, 1)}\phi_t|^2 = o(t)$.
\end{proof}

Claim $\eqref{I_2}$ gives that 
$$I_3(t, x) = -c_1\delta(x)^{d + 2 - n} + c_1c_2^{-1}\delta(x)^{d + 2 + \alpha - n}\ilim_{\R^d}{\frac{\left(t\frac{\partial \phi_t}{\partial t}(0, y) - \phi_t(y)\right)dy}{|x - y|^{d + \alpha}}} + o(t).$$
We also have that
$$I_1(t, x) = c_1\delta(x)^{d + 2 - n} + \ilim_{\R^d}{\frac{\left(\phi_t(y) - t\frac{\partial \phi_t}{\partial t}(0, y)\right)dy}{|x - y|^{n - 2}}}.$$
Conditions 2 (the definition of a Frechet differential) and 4 in Definition $\ref{defn1}$, combined with $\eqref{f_bmo4}$ for $m = 1$, give that
$$c_1c_2^{-1}\ilim_{\R^d}{\frac{\left(t\frac{\partial \phi_t}{\partial t}(0, y) - \phi_t(y)\right)dy}{|x - y|^{d + \alpha}}} \; \mbox{and} \; \ilim_{\R^d}{\frac{\left(\phi_t(y) - t\frac{\partial \phi_t}{\partial t}(0, y)\right)dy}{|x - y|^{n - 2}}} = o(t),$$
and therefore
$$\left|I(t, x) - t \frac{\partial I(t, x)}{\partial t}(0)\right| = \left|I_1(t, x) + I_3(t, x)\right| = o(t).$$
\end{proof}

So we finally proved that $\eqref{take_diff}$ is differentiable at zero, and by $\eqref{t'}$ the derivative is 
$$
    \frac{\partial I(t, x)}{\partial t}(0) = \ilim_{\R^d}{\frac{\frac{\partial \phi_t}{\partial t}(0, y)dy}{|x - y|^{n - 2}}} - c_1c_2^{-1}\delta(x)^{d = 2 + \alpha - n}\ilim_{\R^d}{\frac{\frac{\partial \phi_t}{\partial t}(0, y)dy}{|x - y|^{d + \alpha}}}.
$$
Since $I(t, x)$ is identically zero, the derivative also has to be zero, and therefore, by definition of the constant $c_2$, we have the equation
\begin{equation}\label{flat}
    \ilim_{\R^d}{\left[c_2^{-1}c_1\frac{\delta(x)^{d + 2 + \alpha - n}}{|x - y|^{d + \alpha}} - \frac{1}{|x - y|^{n - 2}}\right]\frac{\partial \phi_t}{\partial t}(0, y)dy}  = 0 \quad \forall \, x \in \R^n \setminus \R^d.
\end{equation}

We claim that $\eqref{flat}$ implies the derivative $\frac{\partial \phi_t}{\partial t}(\cdot, 0)$ has to be a constant function, which is a contradiction with the assumption 3 in Definition $\ref{defn1}$ of non-trivial one-parameter family. Denote $\phi(y) = \frac{\partial \phi_t}{\partial t}(y, 0)$ for simplicity.

\begin{lm}\label{flat3_lm}
Suppose that $\phi(y)$ is a function in $BMO(\R^d)$. Then the equation
\begin{equation}\label{flat2}
    \ilim_{\R^d}{\left[c_2^{-1}c_1\frac{\delta(x)^{d + 2 + \alpha - n}}{|x - y|^{d + \alpha}} - \frac{1}{|x - y|^{n - 2}}\right]\phi(y)dy}  = 0 \quad \forall \, x \in \R^n \setminus \R^d.
\end{equation}
implies that $\phi$ is a constant function.
\end{lm}
\begin{proof}
The idea of the proof is that we can interpret the equation $\eqref{flat2}$ as vanishing of the functional $\phi$ on a certain space of functions. Denote
$$Ker(y) = \frac{c_2^{-1}c_1}{(1 + |y|^2)^{(d + \alpha)/2}} - \frac{1}{(1 + |y|^2)^{(n - 2)/2}}, \; y \in \R^d.$$
Note that $\ilim_{\R^d}{Ker(y)dy} = 0$. Consider the (closed) functional space $\mathcal{S}$ generated by linear combinations with coefficients ${c_i}$ of $L^1$-normalized translations and dilatations of the function $Ker$ with the norm $\|\cdot\| = \inf{\slim_i{|c_i|}}$, where $\inf$ is taken over all representations of an element of $\mathcal{S}$. If $\phi$ was an element of the dual space of $\mathcal{S}$, then $\eqref{flat2}$ would mean that it acts like a zero functional on the space $\mathcal{S}$. Indeed, let $x = (x_0, r, 0 \dots 0)$ be a point in $\R^n \setminus \R^d$, where $x_0$ is the projection of $x$ to the hyperplane $\R^d$ and $r = \delta(x)$. Observe that
$$c_2^{-1}c_1\frac{\delta(x)^{d + 2 + \alpha - n}}{|x - y|^{d + \alpha}} - \frac{1}{|x - y|^{n - 2}} = r^{2 - n}\left(\frac{c_2^{-1}c_1}{(1 + |y - x_0|^2/r^2)^{(d + \alpha)/2}} - \frac{1}{(1 + |y - x_0|^2/r^2)^{(n - 2)/2}}\right)$$ 
$$ =r^{2 - n} Ker\left(\frac{y - x_0}{r}\right).$$
So, $\eqref{flat2}$ implies that $\phi$ vanishes on all the translations and dilatations of the function $Ker$, and therefore on the whole $\mathcal{S}$. 

Therefore we only need to check that $BMO$ is contained in the space dual to the space $\mathcal{S}$. It is known that the latter space is the homogeneous Besov space $\dot B^0_{1, 1}$, see $\cite{YM}$. Its dual space is the homogeneous Besov space $\dot B^0_{\infty, \infty}$, which is larger than $BMO$. 
\end{proof}

\begin{rem}\label{fourier1}
Note that, if we knew in addition that the Fourier transform of $\phi$ is a well-defined function, then we could apply another, a rather elegant, argument to get the conclusion of Lemma $\ref{flat3_lm}$. 

Indeed, denote $x = (z, h)$, $z = (z_1, \dots, z_d)$, $h = (h_{d + 1}, \dots, h_n)$ and
$$f_h(z - y) = c_2^{-1}c_1\frac{\delta(x)^{d + 2 + \alpha - n}}{|x - y|^{d + \alpha}} - \frac{1}{|x - y|^{n - 2}} = \frac{c_2^{-1}c_1|h|^{d + \alpha + 2 - n}}{(|h|^2 + |z - y|^2)^{(d + \alpha)/2}} - \frac{1}{(|h|^2 + |z - y|^2)^{(n - 2)/2}}.$$
With this notation $\eqref{flat2}$ transforms into a convolutional equation
$$f_h * \phi(z) = 0 \quad \forall z \in \R^d.$$
If the Fourier transform is defined for $\phi$, then the convolutional equation above transforms into
$$\widehat{f_h * \phi} = \hat{f_h}\hat{\phi} = 0.$$
This implies that $\hat{\phi}$ can be nonzero only on the set of the common zeros of functions $\hat{f_h}$ for every $h$ (or, rather, every value of $|h|$). What can this set be? First, observe that $\zeta = 0$ is a zero of $\hat{f_h}$ for every $h$. Second, since $\hat{f_h}$ is radially-symmetric, the set of zeros of $\hat{f_{h_1}}$ and of $\hat{f_{h_2}}$ for $h_1$ and $h_2$ such that $|h_1| = |h_2|$ is the same. Third, suppose that $\zeta$ is a zero for $\hat{f_1}$, then, since $f_h(y) = f_1(y/|h|)$ and $\hat{f_h}(\zeta) = |h|\hat{f_1}(|h|\zeta)$, by the definition of $f_h$, $|h|\zeta$ is a zero for $\hat{f_h}$. The latter observation means that if $\zeta \neq 0$ was a common zero for $\hat{f_h}$ for all $h$, then the set $\R \zeta$ would be also contained in the set of common zeros for $\hat{f_h}$. Thus $\hat{f_h}$ would have to be zero everywhere on $\R^d$ for every $h$, which is clearly false.

We conclude then that $\hat{\phi} = 0$ everywhere on $\R^d$ except for zero. This implies $\hat{\phi} = c\delta_0$ and $\phi$ equal to a constant.
\end{rem}

\subsection{No one-parameter families of Lipschitz graph solutions}

In this subsection we prove Theorem $\ref{mainthm_graph}$. Recall what a non-trivial one-parameter differentiable family of smooth graph perturbations of the flat solution is: see Definition $\ref{defn2}$. Our scheme of proof will be essentially the same as in Subsection 4.1, but the expressions we are forced to work with are more bulky.  

Our initial assumptions are that the functions $\psi_t$ are $Ct$-Lipschitz on $\R^d$. However, it is well-known that if a function $f$ is $\Lambda$-Lipschitz, then it is in $W^{1, \infty}(\R^d)$ and one has $\|\nabla f\|_\infty \le \Lambda$. The Rademacher-Calder\'on theorem gives then that the function $f$ is truly differentiable almost everywhere. So without loss of generality we can assume that every $\psi_t$ is differentiable and moreover that $\|D \psi_t\|_\infty \le Ct$. The latter implies the conditions
\begin{equation}\label{graph_cond}
    \|\frac{\partial \psi_t(y)}{\partial y_i}\|_{BMO} \le Ct \quad \mbox{and} \quad \Big|\ilim_{B(0, 1)}{\frac{\partial \psi_t(y)}{\partial y_i}dy}\Big| \le Ct
 \end{equation}
for all $i$, which we will use later on. Recall that $\eqref{graph_cond}$ means that norms of every component of vectors $\frac{\partial \psi_t(y)}{\partial y_i}$ or moduli of every component of vectors $\ilim_{B(0, 1)}{\frac{\partial \psi_t(y)}{\partial y_i}dy}$ admit the estimate $\eqref{graph_cond}$.

Denote by $\eta_t$ the map $\Id + \psi_t$. According to the formula for the Lebesgue measure on a smooth surface, one can rewrite $\eqref{densities}$ as
\begin{equation}\label{diff_gr}
    J(t, x) = c_3 J_1(t, x) - J_2(t, x) = 0 \quad \forall \; x \in \R^n \setminus Im(\eta_t(\R^d)), \quad \mbox{where}
\end{equation}
$$J_1(t, x) = \left(\ilim_{\R^d}{\frac{(1 + \phi_t(y))^\frac{\alpha}{n - d - 2}\sqrt{|\mydet((D\eta_t)^T D\eta_t)|}dy}{|x - \eta_t(y)^{d + \alpha}|}}\right)^\frac{n - d - 2}{\alpha} \quad \mbox{and}$$
$$J_2(t, x) = \ilim_{\R^d}{\frac{(1 + \phi_t(y))\sqrt{|\mydet((D\eta_t)^T D\eta_t)|}dy}{|x - \eta_t(y)|^{n - 2}}}.$$
We keep to the assumption $\|\phi_t\|_\infty \le 1/2$. Again, we want to take the Frechet derivative of the left-hand side of $\eqref{diff_gr}$ at zero, and the natural candidate for it is
$$\frac{\partial J(t, x)}{\partial t}(0) = c_1c_2^{-1}\frac{n - d - 2}{\alpha}\delta(x)^{\alpha + d + 2 - n}\left[\ilim_{\R^d}{\frac{\frac{\alpha}{n - d - 2}\frac{\partial \phi_t}{\partial t}(0, y)}{|x - y|^{d + \alpha}}dy} - \ilim_{\R^d}{\frac{(d + \alpha)\langle x, \frac{\partial \psi_t}{\partial t}(0, y)\rangle}{|x - y|^{d + \alpha + 2}}dy}\right] $$
\begin{equation}\label{diff_graph}
- \ilim_{\R^d}{\frac{\frac{\partial \phi_t}{\partial t}(0, y)}{|x - y|^{n - 2}}dy} + \ilim_{\R^d}{(n - 2)\frac{\langle x, \frac{\partial \psi_t}{\partial t}(0, y)\rangle}{|x - y|^n}dy}.
\end{equation}
Here we just differentiated formally $J(t, x)$ and used the condition 1 in Definition $\ref{defn2}$: as we shall see, we are lucky that the square roots in $\eqref{diff_gr}$ only give terms of order two. Once again, until the end of the proof of the following lemma, $\frac{\partial J(t, x)}{\partial t}(0)$ is just a notation.
\begin{lm}
\begin{equation}\label{o(t)}
    \left|J(t, x) - t \frac{\partial J(t, x)}{\partial t}(0)\right| = o(t), \quad t \to 0.
\end{equation}
\end{lm}
\begin{proof}
We start with proving a technical fact analogous to $\eqref{I_2}$. The proof is ideologically the same, but looks bulkier, because typically we need to multiply three expansions from Taylor's theorem. Plus, we need to deal with the determinant of metric tensor induced by the map $\eta_t$.
\begin{cl}
$$J_2(t, x) = c_1\delta(x)^{d + 2 - n} + \ilim_{\R^d}{\frac{\phi_t(y)dy}{|x - y|^{n - 2}}} - (n - 2)\ilim_{\R^d}{\frac{\langle x, \psi_t(y)\rangle dy}{|x - y|^n}} + o(t), \quad t \to 0, \quad \mbox{and}$$
$$J_1(t, x) = c_2^\frac{n - d - 2}{\alpha}\delta(x)^{d + 2 - n} + c_2^{\frac{n - d - 2}{\alpha} - 1}\delta(x)^{d + 2 + \alpha - n}\frac{n - d - 2}{\alpha}\left[\ilim_{\R^d}{\frac{\frac{\alpha}{n - d - 2}\phi_t(y)}{|x - y|^{d + \alpha}}}\right.$$
$$\left. - \ilim_{\R^d}{\frac{(d + \alpha)\langle x, \psi_t(y)\rangle}{|x - y|^{d + \alpha + 2}}}\right] + o(t), \quad t \to 0.$$
\end{cl}
\begin{proof}
Let us deal first with the square root of the modulus of the determinant of metric tensor induced by $\eta_t$ in $\eqref{diff_gr}$. Denote by $(\psi_t)_j, j = d + 1, \dots, n$, the projection of $\psi_t$ on the coordinate axis $j$. Then the differential $D \eta_t$ is the linear map $\R^d \to \R^n$ represented by the matrix $n \times d$, where in the first $d$ rows the only non-zero elements are ones on the diagonal:
$$
\begin{pmatrix}
1 &  \dots & 0\\
\vdots & \ddots & \vdots\\
0 & \dots & 1\\
\frac{\partial (\psi_t)_{d + 1}}{\partial y_1} & \dots & \frac{\partial (\psi_t)_{d + 1}}{\partial y_d}\\
\vdots & \vdots & \vdots\\
\frac{\partial (\psi_t)_{n}}{\partial y_1} & \dots & \frac{\partial (\psi_t)_{n}}{\partial y_d}
\end{pmatrix}.
$$

Therefore an element $a_{ij}$ of the $d \times d$ matrix $A = (a_{ij})= (D\eta_t)^T D\eta_t$ is equal to $\chi_{i = j} + \slim_{k = d + 1}^n{\frac{\partial (\psi_t)_k}{\partial y_i}\frac{\partial (\psi_t)_k}{\partial y_j}}$, where $\chi_{i = j}$ is one if $i = j$ and zero otherwise. So, the determinant $\mydet((D\eta_t)^T D\eta_t)$ is of the form
$$1 + M_t(y) + R_t(y),$$
where $M_t(y) = \slim_{i = 1}^n{\slim_{k = d + 1}^n{\left(\frac{\partial (\psi_t)_k}{\partial y_i}\right)^2}}$ and
$R_t(y)$ is the rest of the determinant: the sum of products which consist of more than $2(n - d)$ partial derivatives. The point in decomposing the determinant in such a way is that the term $M_t(y)$ has the $BMO$ norm and the average over a unit ball of order $t^{2(n - d)}$, while the $BMO$ norm and the average over the unit ball of the other term $R_t$ decay even faster.

~\

With this at hand, observe that
$$(1 + \phi_t(y))^\frac{\alpha}{n - d - 2} = 1 + \frac{\alpha}{n - d - 2}\phi_t(y) + \phi_t(y)^2 g_1(t, y),$$ 
$$ |\mydet((D\eta_t)^T D\eta_t)|^{1/2} = 1 + (M_t(y) + R_t(y))g_2(t, y),$$
$$|x - \eta_t(y)|^{-(n - 2)} = |x - y|^{-(n - 2)} - (n - 2)|x - y|^{-n}\langle x, \psi_t(y)\rangle + (\langle x, \psi_t(y)\rangle^2 + |\psi_t(y)|^2) |x - y|^{-n} g_3(t, y), \quad \mbox{and}$$
$$|x - \eta_t(y)|^{-(d + \alpha)} = |x - y|^{-(d + \alpha)} - (d + \alpha)|x - y|^{-(d + \alpha + 2)}\langle x, \psi_t(y)\rangle$$ 
$$+ (\langle x, \psi_t(y)\rangle^2 + |\psi_t(y)|^2) |x - y|^{-(d + \alpha + 2)}g_4(t, y),$$
where $g_i(t, y)$ are uniformly bounded in $t$ functions. 

~\

We deal with $J_2(t, x)$ first. The product under the first integral in $J_2(t, x)$ 
$$(1 + \phi_t(y))(1 + (M_t(y) + R_t(y))g_2(t, y))(|x - y|^{-(n - 2)}$$ 
$$ - (n - 2)|x - y|^{-n}\langle x, \psi_t(y)\rangle + (|\psi_t|^2 + \langle x, \psi_t(y)\rangle^2) |x - y|^{-n)}g_3(t, y))$$
is equal to 
$$|x - y|^{-(n - 2)} + \phi_t(y)|x - y|^{-(n - 2)} - (n - 2)|x - y|^{-n}\langle x, \psi_t(y)\rangle - \phi_t(y)(n - 2)|x - y|^{-n}\langle x, \psi_t(y)\rangle$$ 
$$
+ (\psi_t(y)^2 + \langle x, \psi_t(y)\rangle^2 + (M_t(y) + R_t(y)))|x - y|^{-(n - 2)}G_t(y),
$$
where $G_t(y)$ is a uniformly bounded in $t$ function. Integrating, we have
$$J_2(t, x) = c_1\delta(x)^{d + 2 - n} + \ilim_{\R^d}{\frac{\phi_t(y)dy}{|x - y|^{n - 2}}} - (n - 2)\ilim_{\R^d}{\frac{\langle x, \psi_t(y)\rangle dy}{|x - y|^n}}$$
$$ - (n - 2)\ilim_{\R^d}{\frac{\phi_t(y)\langle x, \psi_t(y)\rangle dy}{|x - y|^n}} + \ilim_{\R^d}{\frac{(\psi_t(y)^2 + \langle x, \psi_t(y)\rangle^2 + (M_t(y) + R_t(y)))G_t(y)dy}{|x - y|^{n - 2}}}.$$
By $\eqref{f_bmo4}$ and $\eqref{Ho}$, combined with condition 4 in Definition $\ref{defn2}$ and $\eqref{graph_cond}$, the last two terms are $o(t)$.

~\

Now we deal with $J_1(t, x)$. The product under the integral in $J_1(t, x)$
$$\left(1 + \frac{\alpha}{n - d - 2}\phi_t(y) + \phi_t(y)^2g_1(t, y)\right)\left(1 + (M_t(y) + R_t(y))g_2(t, y))(|x - y|^{-(d + \alpha)}\right.$$ 
$$ \left. - (d + \alpha)|x - y|^{- (d + \alpha + 2)}\langle x, \psi_t(y)\rangle + (\langle x, \psi_t(y)\rangle^2 + |\psi_t(y)|^2) |x - y|^{-(d + \alpha + 2)}g_4(t, y)\right)$$
is equal to
$$|x - y|^{-(d + \alpha)} + \frac{\alpha}{n - d - 2}\phi_t(y)|x - y|^{- (d + \alpha)} - (d + \alpha)|x - y|^{- (d + \alpha + 2)}\langle x, \psi_t(y)\rangle $$
$$
- \frac{\alpha}{n - d - 2}\phi_t(y)(d + \alpha)|x - y|^{- (d + \alpha + 2)}\langle x, \psi_t(y)\rangle$$ $$ + (\psi_t(y)^2 + \langle x, \psi_t(y)\rangle^2 + (M_t(y) + R_t(y)))|x - y|^{-(d + \alpha)}H_t(y),
$$
where $H_t(y)$ is a uniformly bounded in $t$ function. Therefore we have
$$J_1(t, x) = c_2^\frac{n - d - 2}{\alpha}\delta(x)^{d + 2 - n} + c_2^{\frac{n - d - 2}{\alpha} - 1}\delta(x)^{d + 2 + \alpha - n}\ilim_{\R^d}{\frac{\phi_t(y)dy}{|x - y|^{d + \alpha}}}$$ 
$$- c_2^{\frac{n - d - 2}{\alpha} - 1}\delta(x)^{d + 2 + \alpha - n}\frac{n - d - 2}{\alpha}\ilim_{\R^d}{\frac{(d + \alpha)\langle x, \psi_t(y)\rangle dy}{|x - y|^{d + \alpha + 2}}}$$ 
$$ + \ilim_{R^d}{\frac{(\phi_t(y)\langle x, \psi_t(y)\rangle + \langle x, \psi_t(y)\rangle^2 + \psi_t(y)^2 + (M_t(y) + R_t(y)))H'_t(y)dy}{|x - y|^{d + \alpha}}},$$
where $H'_t(y)$ is another uniformly bounded in $t$ function. By $\eqref{f_bmo4}$ and $\eqref{Ho}$, combined with condition 4 in Definition $\ref{defn2}$ and $\eqref{graph_cond}$, the last term is $o(t)$. 
\end{proof}

Denote
$$J_3(t, x) = c_3J_1(t, x) - c_1\delta(x)^{d + 2 - n} - c_1c_2^{-1}\frac{n - d - 2}{\alpha}\delta(x)^{\alpha + d + 2 - n}t \cdot $$ $$\cdot\left[\ilim_{\R^d}{\frac{\frac{\alpha}{n - d - 2}\frac{\partial \phi_t}{\partial t}(0, y)}{|x - y|^{d + \alpha}}dy} - \ilim_{\R^d}{\frac{(d + \alpha)\langle x, \frac{\partial \psi_t}{\partial t}(0, y)\rangle}{|x - y|^{d + \alpha + 2}}dy}\right]$$
and
$$J_4(t, x) = J_2(t, x) - c_1\delta(x)^{d + 2 - n} - t\ilim_{\R^d}{\frac{\frac{\partial \phi_t}{\partial t}(0, y)}{|x - y|^{n - 2}}dy} + t\ilim_{\R^d}{(n - 2)\frac{\langle x, \frac{\partial \psi_t}{\partial t}(0, y)\rangle}{|x - y|^n}dy}.$$
In this notation we have 
$$\left|J(t, x) - t \frac{\partial J(t, x)}{\partial t}(0)\right| = |J_3(t, x) - J_4(t, x)|,$$
because the constant terms cancel out. It is only left to show that the claim above gives $|J_{k}(t, x)| = o(t)$ for $k = 3, 4$, which finishes the proof. Indeed, the claim asserts that
$$J_4(t, x) = \ilim_{\R^d}{\frac{\phi_t(y)dy}{|x - y|^{n - 2}}} - (n - 2)\ilim_{\R^d}{\frac{\langle x, \psi_t(y)\rangle dy}{|x - y|^n}} + o(t) - $$
$$- t\ilim_{\R^d}{\frac{\frac{\partial \phi_t}{\partial t}(0, y)}{|x - y|^{n - 2}}dy} + t\ilim_{\R^d}{(n - 2)\frac{\langle x, \frac{\partial \psi_t}{\partial t}(0, y)\rangle}{|x - y|^n}dy}.$$
But the condition 2 in Definition $\ref{defn2}$, combined with $\eqref{f_bmo4}$ applied with $m = 1$, gives that
$$\ilim_{\R^d}{\frac{\phi_t(y)dy}{|x - y|^{n - 2}}} - t\ilim_{\R^d}{\frac{\frac{\partial \phi_t}{\partial t}(0, y)}{|x - y|^{n - 2}}dy} = o(t) \quad \mbox{and}$$
$$\ilim_{\R^d}{\frac{\langle x, \psi_t(y)\rangle dy}{|x - y|^n}} - t\ilim_{\R^d}{\frac{\langle x, \frac{\partial \psi_t}{\partial t}(0, y)\rangle}{|x - y|^n}dy} = o(t).$$
The integral $J_3(t, x)$ can be treated the same way.
\end{proof}

~\ 

We continue studying $\eqref{diff_gr}$. Since $J(t, x)$ is identically zero, the derivative $\frac{\partial J(t, x)}{\partial t}$, given by $\eqref{diff_graph}$ at zero vanishes as well, which gives
$$\ilim_{\R^d}{\frac{\partial \phi}{\partial t}(0, y)\left[\frac{1}{|x - y|^{n - 2}} - \frac{c_1c_2^{-1}\delta(x)^{\alpha + d + 2 - n}}{|x - y|^{d + \alpha}}\right]dy} $$
$$    = \ilim_{\R^d}{\langle \frac{\partial \psi_t}{\partial t}(0, y), x \rangle\left[\frac{n - 2}{|x - y|^n} - \frac{c_1c_2^{-1}(d + \alpha)\frac{n - d - 2}{\alpha}\delta(x)^{\alpha + d + 2 - n}}{|x - y|^{d + \alpha + 2}}\right]dy}.$$
Without loss of generality we can assume that the projection $\frac{\partial \psi_t}{\partial t}_{d + 1}(0, y)$ of $\frac{\partial \psi_t}{\partial t}(0, y)$ to the axis $d + 1$ is not a constant function. Assume that $x$ is of the form $(x_1, \dots, x_d, x_{d + 1}, 0, \dots)$. Then, with the notation $\psi(y) =  \frac{\partial \psi_t}{\partial t}_{d + 1}(0, y)$ and $\phi(y) = \frac{\partial \phi}{\partial t}(0, y)$, the equation above transforms into
$$\ilim_{\R^d}{\phi(y)\left[\frac{1}{|x - y|^{n - 2}} - \frac{c_1c_2^{-1}\delta(x)^{\alpha + d + 2 - n}}{|x - y|^{d + \alpha}}\right]dy} = $$
\begin{equation}\label{homo_ker}
    \ilim_{\R^d}{\psi(y)x_{d + 1}\left[\frac{n - 2}{|x - y|^n} - \frac{c_1c_2^{-1}(d + \alpha)\frac{n - d - 2}{\alpha}\delta(x)^{\alpha + d + 2 - n}}{|x - y|^{d + \alpha + 2}}\right]dy} \quad \forall x \in \R^{d + 1}. 
\end{equation}
\begin{lm}\label{blm2}
The equation $\eqref{homo_ker}$ implies that the function $\psi$ is zero.
\end{lm}
\begin{proof}
With the notation $h = |x_{d + 1}|$ and $x_0 = (x_1, \dots, x_d)$ we can rewrite the equation $\eqref{homo_ker}$ as
$$\ilim_{\R^d}{h^{2 - n}\phi(y)\left[\frac{1}{(1 + |x_0 - y|^2/h^2)^{(n - 2)/2}} - \frac{c_1c_2^{-1}}{(1 + |x_0 - y|^2/h^2)^{(d + \alpha)/2}}\right]dy} = $$
\begin{equation}\label{homo_ker2}
    \ilim_{\R^d}{h^{1 - n}\psi(y)\left[\frac{n - 2}{(1 + |x_0 - y|^2/h^2)^{n/2}} - \frac{c_1c_2^{-1}(d + \alpha)\frac{n - d - 2}{\alpha}}{(1 + |x_0 - y|^2/h^2)^{(d + \alpha + 2)/2}}\right]dy} \quad \forall h > 0, \; x_0 \in \R^d.
\end{equation}
Denote
$$g_h(y) = \frac{1}{(1 + |y|^2/h^2)^{(n - 2)/2}} - \frac{c_1c_2^{-1}}{(1 + |y|^2/h^2)^{(d + \alpha)/2}} \quad \mbox{and}$$
$$f_h(y) = \frac{n - 2}{(1 + |y|^2/h^2)^{n/2}} - \frac{c_1c_2^{-1}(d + \alpha)\frac{n - d - 2}{\alpha}}{(1 + |y|^2/h^2)^{(d + \alpha + 2)/2}}.$$
Since the Fourier transform is defined for all of the terms in $\eqref{homo_ker2}$, we can rewrite the equation on the Fourier-transform side and get that
\begin{equation}\label{graph_fourier}
   h\hat{\phi}(\zeta)\hat{g_h}(\zeta) = \hat{\psi}(\zeta)\hat{f_h}(\zeta) \quad \forall \zeta \in \R^d \; \forall h > 0. 
\end{equation}
\begin{rem}
It is important for our argument that every term in the products in $\eqref{graph_fourier}$ is a well-defined function. This is the reason why in Definition $\ref{defn2}$ we asked the derivatives $\phi$ and $\psi$ to be in $L^1$. But the operators of convolution with the kernels $g_h$ and $f_h$ are very good: for sure they are classical Calder\`on-Zygmund operators (one can check easily that the kernels satisfy the H\"ormander condition). So one could ask some other reasonable regularity from $\phi$ and $\psi$, as long as their Fourier transforms are well-defined functions. So, instead of writing $\phi, \psi \in L^1 \cap BMO(\R^d)$ in Definition $\ref{defn2}$ we could write $\phi, \psi \in L^p \cap BMO(\R^d)$ for any $p \in [1, 2]$: see \cite{LL}, Sections 5.6 and 5.7.
\end{rem}
The rest of the argument relies on asymptotics of the functions $\hat{g_h}$ and $\hat{f_h}$ at $h \approx 0$. First, observe that $\hat{g_h}(\zeta) = \hat{g_1}(h\zeta)$ and $\hat{f_h}(\zeta) = \hat{f_1}(h\zeta)$. Next, the Fourier transforms $\hat{g_h}(\zeta)$ and $\hat{f_h}(\zeta)$ are rather easy to compute. Indeed, one has, by the integral definition of the modified Bessel function of the second kind $K_b(z)$, 
\begin{equation}\label{ft}
    \left(\frac{1}{(1 + |y|^2)^{a/2}}\right)^{\widehat{}}(h\zeta) = \frac{2^{1 - a/2}h^{(a - 1)/2}|\zeta|^{(a - 1)/2}K_{\frac{a - 1}{2}}(h|\zeta|)}{\Gamma\left(\frac{a}{2}\right)},
\end{equation}
where
$$K_b(z) = \frac{\Gamma(b + \frac{1}{2})(2z)^b}{\sqrt{\pi}}\ilim_0^\infty{\frac{\cos(t)dt}{(t^2 + z^2)^{(b + 1)/2}}}.$$
From the series representation of $K_b$ (see, for example, $\cite{AS}$) it follows that, if $b > 0$, for small arguments $z$ one has $K_b(z) \sim \frac{\Gamma(b)}{2}\left(\frac{2}{z}\right)^b$, where by $\sim$ we denote the asymptotic equivalence of functions. Wherefore
$$\left(\frac{1}{(1 + |y|^2)^{a/2}}\right)^{\widehat{}}(h\zeta) \sim \frac{2^{-1/2}\Gamma\left(\frac{a - 1}{2}\right)}{\Gamma\left(\frac{a}{2}\right)}$$
for $\zeta \neq 0$ fixed and $h$ small. Then for fixed $\zeta \neq 0$ and small $h$ we have
$$\hat{g_h}(\zeta) \sim \frac{\Gamma\left(\frac{n - 3}{2}\right)}{\Gamma\left(\frac{n - 2}{2}\right)} - \frac{c_1c_2^{-1}\Gamma\left(\frac{d + \alpha - 1}{2}\right)}{\Gamma\left(\frac{d + \alpha}{2}\right)}, \quad \mbox{and}$$
$$\hat{f_h}(\zeta) \sim  \frac{(n - 2)\Gamma\left(\frac{n - 1}{2}\right)}{\Gamma\left(\frac{n}{2}\right)} - \frac{c_1c_2^{-1}\frac{n - d - 2}{\alpha}(d + \alpha)\Gamma\left(\frac{d + \alpha + 1}{2}\right)}{\Gamma\left(\frac{d + \alpha + 2}{2}\right)}.$$
Since 
$$c_1c_2^{-1} = \frac{\Gamma\left(\frac{n - d - 2}{2}\right)\Gamma\left(\frac{d + \alpha}{2}\right)}{\Gamma\left(\frac{n - 2}{2}\right)\Gamma\left(\frac{\alpha}{2}\right)},$$
$$\hat{g_h}(\zeta) \sim C_g := \Gamma\left(\frac{n - 3}{2}\right)\Gamma\left(\frac{\alpha}{2}\right) - \Gamma\left(\frac{n - d - 2}{2}\right)\Gamma\left(\frac{d + \alpha - 1}{2}\right) \quad \mbox{and}$$
$$\hat{f_h}(\zeta) \sim C_f := \Gamma\left(\frac{n - 1}{2}\right)\Gamma\left(\frac{\alpha + 2}{2}\right) - \Gamma\left(\frac{n - d}{2}\right)\Gamma\left(\frac{d + \alpha + 1}{2}\right).$$
Observe that, for fixed $\alpha > 0$ which is not ``magic'', $C_f$ and $C_g$ are never simultaneously zero except for the case $d = 1$. Indeed, using the relation $\Gamma(z + 1) = z\Gamma(z)$, it is easy to see that
$$C_g = \frac{\Gamma\left(\frac{n - 1}{2}\right)\Gamma\left(\frac{\alpha + 2}{2}\right)}{\frac{n - 3}{2}\frac{\alpha}{2}} - \frac{\Gamma\left(\frac{n - d}{2}\right)\Gamma\left(\frac{d + \alpha + 1}{2}\right)}{\frac{n - d - 2}{2}\frac{d + \alpha - 1}{2}}.$$
Since $\Gamma\left(\frac{n - 1}{2}\right)\Gamma\left(\frac{\alpha + 2}{2}\right) = \Gamma\left(\frac{n - d}{2}\right)\Gamma\left(\frac{d + \alpha + 1}{2}\right)$ when $C_f = 0$, if $C_g$ is also zero, one has 
$(n - d - 2)(d + \alpha - 1) = (n - 3)\alpha$, which is true either for the ``magic'' $\alpha$ or for $d = 1$. 

For the case when $C_g = C_f = 0$ we have to use the next term of the expansion of $K_b$ at $h \approx 0$: $K_b(z) \sim  \frac{\Gamma(b)}{2}\left(\frac{2}{z}\right)^b + \frac{\Gamma(b)}{2(1 - b)}\left(\frac{2}{z}\right)^{b - 2}$. Note that for all the cases when we will use the formula the parameter $b$ in it will be not equal to one. The second term of the expansion gives that, if $C_f = C_g = 0$,
$$\hat{g_h}(\zeta) \sim h^2C_g' \quad \mbox{with} \; C_g' = \frac{\Gamma\left(\frac{n - 3}{2}\right)}{4 - n} - \frac{\Gamma\left(\frac{d + \alpha - 1}{2}\right)\Gamma\left(\frac{n - d - 2}{2}\right)}{(2 - d - \alpha)\Gamma\left(\frac{\alpha}{2}\right)} \quad \mbox{and}$$
$$\hat{f_h}(\zeta) \sim h^2C_f' \quad \mbox{with} \; C_f' = \frac{\Gamma\left(\frac{n - 1}{2}\right)}{2 - n} + \frac{\Gamma\left(\frac{n - d}{2}\right)\Gamma\left(\frac{d + \alpha + 1}{2}\right)}{(d + \alpha)\Gamma\left(\frac{\alpha + 2}{2}\right)},$$
where $C_f'$ and $C_g'$ are not equal to zero if $\alpha \neq n - d - 2$: this is easy to see since $C_f = C_g = 0$ gives $\Gamma\left(\frac{n - 1}{2}\right)\Gamma\left(\frac{\alpha + 2}{2}\right) = \Gamma\left(\frac{n - d}{2}\right)\Gamma\left(\frac{d + \alpha + 1}{2}\right)$ and $\Gamma\left(\frac{n - 3}{2}\right)\Gamma\left(\frac{\alpha}{2}\right) = \Gamma\left(\frac{n - d - 2}{2}\right)\Gamma\left(\frac{d + \alpha - 1}{2}\right)$.

Let us return now to $\eqref{graph_fourier}$ and restrict ourselves to small $h$. Fix $\zeta \neq 0$. Observe that asymptotically the equation $\eqref{graph_fourier}$ looks like
\begin{equation}\label{gr_f_asymp}
    h\hat{\phi}(\zeta)(C_g + h^2C_g') = \hat{\psi}(\zeta)(C_f + h^2C_f') \quad \forall \; h > 0 \; \mbox{small}.
\end{equation}
It is now clear that $\eqref{graph_fourier}$ can never be true unless $\hat{\psi}(\zeta) = 0$. Suppose the contrary. Then we have to have $C_f = 0$ in $\eqref{gr_f_asymp}$. If $\hat{\phi}(\zeta)$ is also not zero, then $C_g, C_g'$ and $C_f'$ are zero in $\eqref{gr_f_asymp}$, but we saw above that this cannot be true. This implies $\hat{\phi}(\zeta) = 0$, but then $C_f'$ has to be zero as well. At the same time, from the computations above, we saw that $C_f$ and $C_f'$ can be simultaneously zero only if $\alpha$ is ``magic''. A contradiction.

So we have that $\hat{\psi}(\zeta) = 0$ everywhere except probably for $\zeta = 0$. Therefore the function $\psi$ can only be a constant, which also contradicts the definition of one-parameter differentiable family of graph perturbations.

\end{proof}

\subsection{A comment on a sporadic family versus one-parameter family, and BMO versus the Besov space}

We would like to give some more comments on why the plan we described in the beginning of this section for solving the hypothesis related to $\eqref{solution}$ stated in the introduction does not work. Essentially, two more technical steps separate the results of Theorems $\ref{mainthm_flat}$ and $\ref{mainthm_graph}$ from a proof of the hypothesis that there are no solutions of the equation $\eqref{solution}$ in a small neighbourhood of the flat solution. To be more specific, let us discuss just the flat case: the case of perturbations of density of a measure on $\R^d$. Recall that if we suppose the hypothesis to be false, we get a discrete family $\{1 + \phi_t\}$ of non-flat solutions of $\eqref{densities}$, because a solutions should exist arbitrarily close to the flat one.

The first step is breaching the gap between such a family of solutions, which we will call sporadic, and a family of perturbations with a continuous parameter $t$ (recall that the family in Definition $\ref{defn1}$ is parametrized by a segment $[0, t_0)$). We claim that this step would not be a big deal. An attentive look at Definition $\ref{defn1}$ and the proof of Theorem $\ref{mainthm_flat}$ shows that continuity of the parameter $t$ is not really needed anywhere. What we really need is the existence of a non-trivial limit $F$ of the sequence $\{\frac{1}{t}\phi_t\}$, indexed by an arbitrary family of points $t$, not necessarily a segment, -- a surrogate for the Fr\'echet derivative. The Banach-Alaoglu theorem, which we mentioned above already, would guarantee the existence as long as the estimate of the sort $\frac{1}{t}\|\phi_t\| \sim 1$ is provided. The latter would not be a problem: we would just need to say that the value of the parameter $t$ is morally the same as $\|\phi_t\|$, where the norm could be a norm in any space we need.

But we would also need to make the second step and to prove that the limit $F$ is a constant. This, as far as the author can see, seems to be (too) hard to do. In Subsection $5.1$ we are able to do this because we fixed in advance the functional space we work in -- the $BMO$ space. If $\{\frac{1}{t}\phi_t\}$ is ``differentiable'' in $BMO$, we can prove that $\eqref{densities}$ is also differentiable and conclude that, since $\eqref{flat}$ is true for $F$ instead of $\frac{\partial \phi_t}{\partial t}(0, y)$, $F$ is zero in $BMO$. However for an arbitrary sporadic family of perturbations we have no a priori indications that the ``derivative'' $F$ is a $BMO$ function. Moreover, the ``smallest'' space possible for it to fit in is indicated by the equation $\eqref{flat}$. Recall that this equation says that $F$ acts like a zero functional on the space $\mathcal{S}$ generated by linear combinations of the family of functions $\{c_2^{-1}c_1\frac{\delta(x)^{d + 2 + \alpha - n}}{|x - y|^{d + \alpha}} - \frac{1}{|x - y|^{n - 2}}\}_{x \in \R^n \setminus \R^d}$. This would imply that $F$ is a constant if $F$ was an element of the space dual to the space $\mathcal{S}$. So, we would do just fine with the assumption that $F$ is an element of $BMO$ if $S$ was the predual to $BMO$ (or would be smaller). But we saw already that the latter is not true: the space $S$ is the homogeneous Besov space $\dot B^0_{1, 1}$, and its dual is the homogeneous Besov space $\dot B^0_{\infty, \infty}$, which is larger than $BMO$. Therefore the strongest assumption we could have made without loss of generality in our argument is that $F$ lives in $\dot B^0_{\infty, \infty}$.
But for this space we have certain indications that one cannot prove that the existence of the ``derivative'' $F$ implies differentiability of $\eqref{densities}$. Namely, the estimates of the type $\eqref{f_bmo4}$ for the appropriate norm ($\|\cdot\|_{\dot B^0_{\infty, \infty}}$ instead of $\|\cdot\|_{BMO}$) are false, morally, because functions from the Besov space are not locally integrable. 

\section{No smooth solutions for the hyperplane}

In this section we describe an attempt to solve directly the equation $\eqref{densities}$ on the hyperplane $\R^d$ and we prove Theorem $\ref{mainthm_direct}$ from Introduction, which asserts that there are no solutions $D_{\alpha, \mu}$ of $\eqref{solution}$ for $E = \R^d$ among the measures $\mu$ with densities of class $C^{2, \varepsilon}(\R^d)$ for any $0 < \varepsilon < \varepsilon_0$ which are not constants.

Let $D_{\alpha, \mu}$ be a flat solution of the equation $L_\alpha D_\alpha = 0$ where $E$ is $\R^d$, and $f \in L^\infty(\R^d)$ be the density of the measure $\mu$. Then the function $h = f^\frac{n - d - 2}{\alpha}$, as discussed in the beginning of Section 5, satisfies the equation 
\begin{equation}\label{fsolution}
\delta(x)^{n - d - 2}\ilim_{\R^d}{\frac{h(y)dy}{|x - y|^{n - 2}}} = c_3\left(\delta(x)^\alpha\ilim_{\R^d}{\frac{h(y)^\frac{\alpha}{n - d - 2}dy}{|x - y|^{d + \alpha}}}\right)^\frac{n - d - 2}{\alpha}, \quad \forall \; x \in \R^n \setminus \R^d.
\end{equation}
This is just $\eqref{densities}$ multiplied by the correct power of $\delta(x)$.

We assume from now on that $h \in C^{2, \varepsilon}(\R^d)$ for a fixed $\varepsilon$. Then we can write a Taylor expansion at an arbitrary point $y_0 \in \R^d$:
$$h(y) = h(y_0) + \langle \nabla h(y_0), y - y_0\rangle + (y - y_0)^T \Hess h(y_0) (y - y_0) + o(|y - y_0|^2), \quad \mbox{and}$$
$$h(y)^\frac{\alpha}{n - d - 2} = h(y_0)^\frac{\alpha}{n - d - 2} + h(y_0)^{\frac{\alpha}{n - d - 2} - 1}\frac{\alpha}{n - d - 2}\left(\langle \nabla h(y_0), y - y_0\rangle \right.$$ 
$$\left. + (y - y_0)^T \Hess h(y_0) (y - y_0)\right) + h(y_0)^{\frac{\alpha}{n - d - 2} - 2}\frac{\frac{\alpha}{n - d - 2}\left(\frac{\alpha}{n - d - 2} - 1\right)}{2}\langle \nabla h(y_0), y - y_0\rangle^2 + o(|y - y_0|^2).$$
Choose $x = (y_0, r, 0 \dots), y_0 \in \R^d$, $0 < r < 1$, and rewrite the equation $\eqref{fsolution}$, using the expansions from above. On the left-hand side we get
$$c_1h(y_0) + r^{n - d - 2}\Big\langle\nabla h(y_0), \ilim_{\R^d}{\frac{(y - y_0)dy}{|(y_0, r) - y|^{n - 2}}}\Big\rangle $$ 
\begin{equation}\label{lhs6}
  + r^{n - d - 2}\ilim_{\R^d}{\frac{(y - y_0)^T\Hess h(y_0)(y - y_0)dy}{|(y_0, r) - y|^{n - 2}}} + r^{n - d - 2}\ilim_{\R^d}{\frac{o(|y - y_0|^2)dy}{|(y_0, r) - y|^{n - 2}}},  
\end{equation}
and on the right-hand side inside the brackets with the power $\frac{n - d - 2}{\alpha}$ we get
$$c_2h(y_0)^\frac{\alpha}{n - d - 2} + h(y_0)^{\frac{\alpha}{n - d - 2} - 1}\frac{\alpha}{n - d - 2}r^{\alpha}\Big\langle\nabla h(y_0), \ilim_{\R^d}{\frac{(y - y_0)dy}{|(y_0, r) - y|^{d + \alpha}}}\Big\rangle + h(y_0)^{\frac{\alpha}{n - d - 2} - 1}\frac{\alpha}{n - d - 2}r^{\alpha}$$ 
\begin{equation}\label{rhs6}
    \cdot\ilim_{\R^d}{\frac{(y - y_0)^T\Hess h(y_0)(y - y_0) + h(y_0)^{-1}\frac{1}{2}\left(\frac{\alpha}{n - d - 2} - 1\right)\langle\nabla h(y_0), y - y_0\rangle^2 }{|(y_0, r) - y|^{d + \alpha}}dy}
\end{equation}
$$+ r^\alpha\ilim_{\R^d}{\frac{o(|y - y_0|^2)dy}{|(y_0, r) - y|^{d + \alpha}}}.$$

Then we expand the right-hand side of $\eqref{fsolution}$ using $(c + u)^\frac{n - d - 2}{\alpha} = c^\frac{n - d - 2}{\alpha} + c^{\frac{n - d - 2}{\alpha} - 1}u + o(u)$ for small $u$ with $c = c_2h(y_0)^\frac{\alpha}{n - d - 2}$ and $u = \eqref{rhs6} - c_2h(y_0)^\frac{\alpha}{n - d - 2}$. We know already from Subsection 5.1 that, because of the symmetries of the denominator, 
$$\ilim_{\R^d}{\frac{(y - y_0)dy}{|(y_0, r) - y|^{n - 2}}} = \ilim_{\R^d}{\frac{(y - y_0)dy}{|(y_0, r) - y|^{d + \alpha}}} = 0,$$
so the second term with $\nabla h(y_0)$ both in $\eqref{lhs6}$ and $\eqref{rhs6}$ vanishes. Therefore the Taylor expansion of the right-hand side of $\eqref{fsolution}$ simplifies to  
$$c_1h(y_0) + c_1c_2^{-1}r^\alpha\cdot\ilim_{\R^d}{\frac{(y - y_0)^T\Hess h(y_0)(y - y_0) + h(y_0)^{-1}\frac{1}{2}\left(\frac{\alpha}{n - d - 2} - 1\right)\langle\nabla h(y_0), y - y_0\rangle^2 }{|(y_0, r) - y|^{d + \alpha}}dy} $$
\begin{equation}\label{Rhs6}
    + r^\alpha\ilim_{\R^d}{\frac{o(|y - y_0|^2)dy}{|(y_0, r) - y|^{d + \alpha}}}.
\end{equation}

Now we would like to treat $\eqref{lhs6}$ and $\eqref{Rhs6}$ as functions of the parameter $r$ to ``get rid of'' the terms with $o(|y - y_0|^2)$. Observe that the term with integral of a fraction with $O(|y - y_0|^2)$ in the numerator in $\eqref{lhs6}$ is of the order $r^2$, and the same is true for the term with integral of a fraction with $O(|y - y_0|^2)$ in the numerator in $\eqref{Rhs6}$. 

Next, with the assumption $h \in C^{2, \varepsilon}$ we can estimate the residues in $\eqref{lhs6}$ and $\eqref{Rhs6}$ and say that $|o(|y - y_0|^2)| \le C|y - y_0|^{2 + \varepsilon}$, where $C$ depends only on H\"older coefficients of second-order partial derivatives of $h$. Therefore we can bound the modulus of the last term in $\eqref{lhs6}$ the following way:
$$\left|r^{n - d - 2}\ilim_{\R^d}{\frac{o(|y - y_0|^2)dy}{|(y_0, r) - y|^{n - 2}}}\right| \le Cr^{2 + \varepsilon}\ilim_{\R^d}{\frac{|y - y_0|^{2 + \varepsilon}/r^{2 + \varepsilon} }{(1 + |y - y_0|^2/r^2)^{(n - 2)/2}}r^{-d}dy}.$$
Assuming $n - d > 4$, the integral above is bounded and does not depend on $r$. So the $o(|y - y_0|^2)$ term in $\eqref{lhs6}$ is of the order $r^{2 + \varepsilon}$. The same is true about the $o(|y - y_0|^2)$ term in $\eqref{Rhs6}$, since
$$\left|r^\alpha\ilim_{\R^d}{\frac{o(|y - y_0|^2)dy}{|(y_0, r) - y|^{d + \alpha}}}\right| \le Cr^{2 + \varepsilon}\ilim_{\R^d}{\frac{|y - y_0|^{2 + \varepsilon}/r^{2 + \varepsilon} }{(1 + |y - y_0|^2/r^2)^{(d + \alpha)/2}}r^{-d}dy},$$
and the integral above is once again bounded and does not depend on $r$, if we assume that $\alpha > 2 + \varepsilon_0$.
We divide $\eqref{fsolution}$ by $r^2$, let $r$ tend to zero, and get that then the $r^2$-order terms should match. Thus we get the equation
$$r^{n - d - 2}\ilim_{\R^d}{\frac{(y - y_0)^T\Hess h(y_0)(y - y_0)dy}{|(y_0, r) - y|^{n - 2}}} $$ $$ = \frac{c_1}{c_2} r^\alpha \ilim_{\R^d}{\frac{(y - y_0)^T\Hess h(y_0)(y - y_0) + h(y_0)^{-1}\frac{1}{2}\left(\frac{\alpha}{n - d - 2} - 1\right)\langle\nabla h(y_0), y - y_0\rangle^2 }{|(y_0, r) - y|^{d + \alpha}}dy}.$$

We get rid of all second-order derivatives with different indices because of the symmetries, and also of all the products of first-order derivatives with different indices as well. Passing to non-tangential limits will give us the following PDE equation for the function $h$:
\begin{equation}\label{PDE0}
    \tilde{c_1}\Delta h = \frac{c_1}{c_2}\tilde{c_2}\Delta h + \frac{1}{2}\frac{c_1}{c_2}\tilde{c_2}\left(\frac{\alpha}{n - d - 2} - 1\right)h^{-1}|\nabla h|^2 \quad \mbox{on} \; \R^d,
\end{equation}
where 
$$\tilde{c_1} = \Vol(\Sf^{d - 1})\ilim_0^\infty{\frac{x^{d + 1}dx}{(1 + x^2)^\frac{n - 2}{2}}} = \Vol(\Sf^{d - 1})\frac{1}{2}\frac{\Gamma\left(\frac{d + 2}{2}\right)\Gamma\left(\frac{n - d - 4}{2}\right)}{\Gamma\left(\frac{n - 2}{2}\right)}, \quad \mbox{and}$$ 
$$\tilde{c_2} = \Vol(\Sf^{d - 1})\ilim_0^\infty{\frac{x^{d + 1}dx}{(1 + x^2)^\frac{d + \alpha}{2}}} = \Vol(\Sf^{d - 1})\frac{1}{2}\frac{\Gamma\left(\frac{d + 2}{2}\right)\Gamma\left(\frac{\alpha - 2}{2}\right)}{\Gamma\left(\frac{d + \alpha}{2}\right)}.$$
Recall that
$$c_1 = \Vol(\Sf^{d - 1})\frac{1}{2}\frac{\Gamma\left(\frac{d}{2}\right)\Gamma\left(\frac{n - d - 2}{2}\right)}{\Gamma\left(\frac{n - 2}{2}\right)} \quad \mbox{and} \quad c_2 = \Vol(\Sf^{d - 1})\frac{1}{2}\frac{\Gamma\left(\frac{d}{2}\right)\Gamma\left(\frac{\alpha}{2}\right)}{\Gamma\left(\frac{d + \alpha}{2}\right)},$$
so
$$\frac{\tilde{c_1}}{\tilde{c_2}}\frac{c_2}{c_1} = \frac{\Gamma\left(\frac{n - d - 4}{2}\right)\Gamma\left(\frac{\alpha}{2}\right)}{\Gamma\left(\frac{n - d - 2}{2}\right)\Gamma\left(\frac{\alpha - 2}{2}\right)}=\frac{\alpha - 2}{n - d - 4}.$$
Therefore the equation $\eqref{PDE0}$ is trivial when the parameter $\alpha$ is ``magic'' and equal to $n - d - 2$. Otherwise $\eqref{PDE0}$ gives an equation of the form
\begin{equation}\label{PDE}
\Delta h = - C h^{-1}|\nabla h|^2 \quad \mbox{on} \; \R^d,    
\end{equation}
where the constant $C$, as one can see easily from $\eqref{PDE0}$, is
$$C = \frac{1}{2}\left(\frac{\alpha}{n - d - 2} - 1\right)\frac{\frac{c_1}{c_2}\tilde{c_2}}{\frac{c_1}{c_2}\tilde{c_2} - \tilde{c_1}} = \frac{1}{2}\left(\frac{\alpha}{n - d - 2} - 1\right)\frac{1}{1 - \frac{\tilde{c_1}}{\tilde{c_2}}\frac{c_2}{c_1}}.$$
Given the computation above, we have
\begin{equation}\label{const}
C = - \left(\frac{1}{2} - \frac{1}{n - d - 2}\right).
\end{equation}

~\

We claim that we can find a change of variables that transforms the equation $\eqref{PDE}$ into the equation $\Delta \cdot = 0$. Indeed, take $g = h^\beta, \, \beta \neq 0, 1$, then
$$\nabla g = \beta h^{\beta - 1}\nabla h, \quad \Delta g = \beta h^{\beta - 1}\Delta h + \beta(\beta - 1) h^{\beta - 2}|\nabla h|^2.$$
Substituting this in $\eqref{PDE}$, we get that $g$ is a solution of 
$$\Delta g = -\frac{|\nabla g|^2}{g}\left(\frac{C + 1}{\beta} - 1\right).$$
So, if $C \neq -1$, we choose $\beta = C + 1$ and we are done. Otherwise take $g = \log{h}$, then
$$h = e^g, \quad \nabla g = h^{-1}\nabla h, \quad \Delta g = h^{-1}\Delta h - h^{-2}|\nabla h|^2.$$
Substituting this in $\eqref{PDE}$, we get that $g$ is a solution of
$$\Delta g = -(C + 1)|\nabla g|^2,$$
and since $C = -1$, the function $g$ is harmonic as we wanted.

So we have that some power of the function $h$ or $\log{h}$ is a harmonic function $g$ on the whole $\R^d$. Since $h$ is the density of an Ahlfors-regular measure on $\R^d$, $h$ is bounded and bounded away from zero. Thus $g$ is in addition bounded and bounded below. So-called Liouville's theorem implies that $g$ has to be a constant, therefore $h$ is a constant as well.

~\

\begin{rem}
In the end of this section we would like to comment on the regularity restrictions we posed on the density $h$ of a measure $\mu$ in Theorem $\ref{mainthm_direct}$. Asking $h$ to be in $C^{2, \varepsilon}(\R^d)$ seems clumsy. But we cannot significantly release the restrictions, because in the core of the method we use lies the differential equation $\eqref{PDE}$. It looks like one cannot interpret this equation without some boundedness assumptions on its right-hand side. 

One could still argue that, given the method with all the changes of variables we implemented to solve $\eqref{PDE}$, if in $\eqref{fsolution}$ we considered $h^\beta$ with some $\beta$ instead of $h$, probably we could get a simpler equation, which requires less regularity assumptions. Indeed, this is clearly true from the proof given above: there is (for almost all values of the constant $C$) a $\beta$ such that, at the end of the day, we would get that $h^\beta$ is a positive harmonic function, bounded and bounded away from zero. However, to get this equation we still need to use the Taylor expansion and to show that the terms in equations like $\eqref{lhs6}$ and $\eqref{rhs6}$ which correspond to small-o terms in the expansion can be neglected. This seems to be hard to do without the regularity assumption we asked for.
\end{rem}

~\

\noindent Polina Perstneva\\
Universit\'e Paris Saclay, LMO \\
polina.perstneva@universite-paris-saclay.fr

\end{document}